# RATES OF CONVERGENCE OF SOME MULTIVARIATE MARKOV CHAINS WITH POLYNOMIAL EIGENFUNCTIONS

By Kshitij Khare and Hua Zhou

*Stanford University and University of California, Los Angeles*

We provide a sharp nonasymptotic analysis of the rates of convergence for some standard multivariate Markov chains using spectral techniques. All chains under consideration have multivariate orthogonal polynomial as eigenfunctions. Our examples include the Moran model in population genetics and its variants in community ecology, the Dirichlet-multinomial Gibbs sampler, a class of generalized Bernoulli–Laplace processes, a generalized Ehrenfest urn model and the multivariate normal autoregressive process.

**1. Introduction.** The theory of Markov chains is one of the most useful tools of applied probability and has numerous applications. Markov chains are used for modeling physical processes and evolution of a population in population genetics and community ecology. Another important use is simulating from an intractable probability distribution. It is a well-known fact that under mild conditions discussed in [2], a Markov chain converges to its stationary distribution. In the applications mentioned above, often it is useful to know exactly how many steps it takes for a Markov chain to be reasonably close to its stationary distribution. Answering this question as accurately as possible is what finding "rates of convergence" of Markov chains is about.

In the current paper, we provide a sharp nonasymptotic analysis of rates of convergence to stationarity for a variety of multivariate Markov chains. This helps determine exactly what number of steps is necessary and sufficient for convergence. These Markov chains appear as standard models in population genetics, ecology, statistics and image processing.

Here is an example of our results. In community ecology, scientists study diversity and species abundances in ecological communities. The Unified

Neutral Theory of Biodiversity and Biogeography (UNTB) is an important theory proposed by ecologist Stephen Hubbell in his monograph [27]. There are two levels in Hubbell's theory, a metacommunity and a local community.

The metacommunity has constant population size $N_M$ and evolves as follows. At each step, a randomly chosen individual is replaced by a new one. With probability $s$ (speciation), the new individual is a new species that never occurs before. With probability $1-s$ (no speciation), the new individual is a copy of one (randomly chosen) of the remaining $N_M - 1$ individuals. After a sufficiently long time, the metacommunity reaches equilibrium which is the celebrated Ewen's sampling distribution. This process can be considered as a variant of the so-called infinite-allele Moran model in population genetics [18]. In population genetics, the reproducing individual could be the same as the dying one.

The local community has constant population size $N$, which is much smaller than $N_M$ in scale. The evolution of the local community is similar to the metacommunity except it has migration instead of speciation. Specifically, at each step, one individual is randomly chosen to die. With probability $m$ (migration), the new individual is an immigrant randomly chosen from the metacommunity. With probability $1-m$ (no migration), the new individual will be a copy of one (randomly chosen) of the remaining $N-1$ individuals in the local community. Again this process is a variant of the so-called multi-allele Moran model in population genetics [18]. The metacommunity evolves at a much larger time scale and is assumed to be fixed during the evolution of the local community.

Since the publication of [27], UNTB received both acclaims and criticisms. For example, in [38], ecologist McGill tests some classical datasets against the equilibrium species abundance of the local community predicted by the UNTB. McGill raised the "time to equilibrium" issue which is of both practical and scientific interests. In order to generate the equilibrium distribution of the local community, he actually simulates the evolution process of the local community on computer. As he claims ([38], page 882):

> No indication of what fixed number of time steps is appropriate has been published. My own experiments show that it can be large...

Also scientifically it is desirable to know how soon a local community reaches equilibrium. One hundred years? One thousand years?

A simulation of the local community process is performed on computer. Suppose that the metacommunity has $d = 5$ species with uniform species frequencies $\mathbf{p} = (0.2, 0.2, 0.2, 0.2, 0.2)$. The local community has population size $N = 20$ and the migration probability is $m = 0.05$. Assume that initially all 20 individuals in the local community are of the same species. Five thousand independent replicas of the local community process are simulated for 1000 steps. For any (random) count vector $\mathbf{X} = (X_1, \ldots, X_d)$ of the local community, where $X_i$ is the count of individuals of species $i$, we can define a Watterson type statistics (population homogeneity) as

$$W(\mathbf{X}) = \sum_{i=1}^{d} \frac{X_i^2}{N^2}.$$

The empirical distributions of the Watterson statistics are plotted along time in Figure 1. By visual inspection, we *suspect* that the distribution of Watterson statistics is close to its stationary distribution after a few hundred steps. A commonly used measure for distance to stationarity for Markov chains is the chi-square distance. Consider a Markov chain with state space $\mathcal{X}$, transition density $K(\cdot,\cdot)$ and stationary density $m$ with respect to a $\sigma$-finite measure $\mu$. The chi-square distance to stationarity after $\ell$ steps, starting at state $x$, is defined as

$$\chi_x^2(\ell) = \int_{\mathcal{X}} \frac{[K^\ell(x,x') - m(x')]^2}{m(x')} \mu(dx').$$

Proposition 4.10 provides rigorous quantitative answers to the question of how long it takes the local community process to be close to its stationary distribution. Let $\chi^2_{N\mathbf{e}_i}(l)$ be the chi-square distance between the distribution of the $d$-dimensional vector of local community species abundances after $\ell$ steps and the stationary distribution, assuming initially all individuals are of species $i$. Proposition 4.10 tells us that 595 steps are necessary and sufficient for convergence. This means that, for $l \leq 595$ steps, $\chi^2_{N\mathbf{e}_i}(l)$ is high and, for $l \geq 595 + 190c$ ($c$ is any positive constant) steps, $\chi^2_{N\mathbf{e}_i}(l) \leq e^{-c}$. For example, by $l = 595 + 190 \times \log 100 \approx 1470$ steps, the chi-square distance $\chi^2_{N\mathbf{e}_i}(l) \leq 0.01$. If this is a tree population with mortality rate 1% per year, 595 and 1470 steps translate into 2975 and 7350 years, respectively. For people who prefer total variation distance, it should be kept in mind that the chi-square distance always produces an upper bound for the total variation distance (see Section 2.1). Figure 2 shows how the (exact) chi-square distance for the $d$-dimensional Markov chain is decreasing over time.

The calculations work because the Markov chain corresponding to the local community process admits a system of multivariate orthogonal polynomials as eigenfunctions. Then a summation formula due to Griffiths (reviewed in Section 2.2.2) pertinent to this system of orthogonal polynomials allows us to do explicit calculations for the chi-square distance.

We provide similar results for all other examples considered in this paper. For every Markov chain, we find positive constants $D$ and $R$, which depend on various parameters of the Markov chain, such that after $D - cR$ steps the chi-square distance is larger than an explicit constant multiple of $e^c$, and after $D + cR$ steps the chi-square distance to stationarity is less than an explicit constant multiple of $e^{-c}$ or similar simple functions. In this sense, we say that $D$ steps are necessary and sufficient for convergence in chi-square distance.

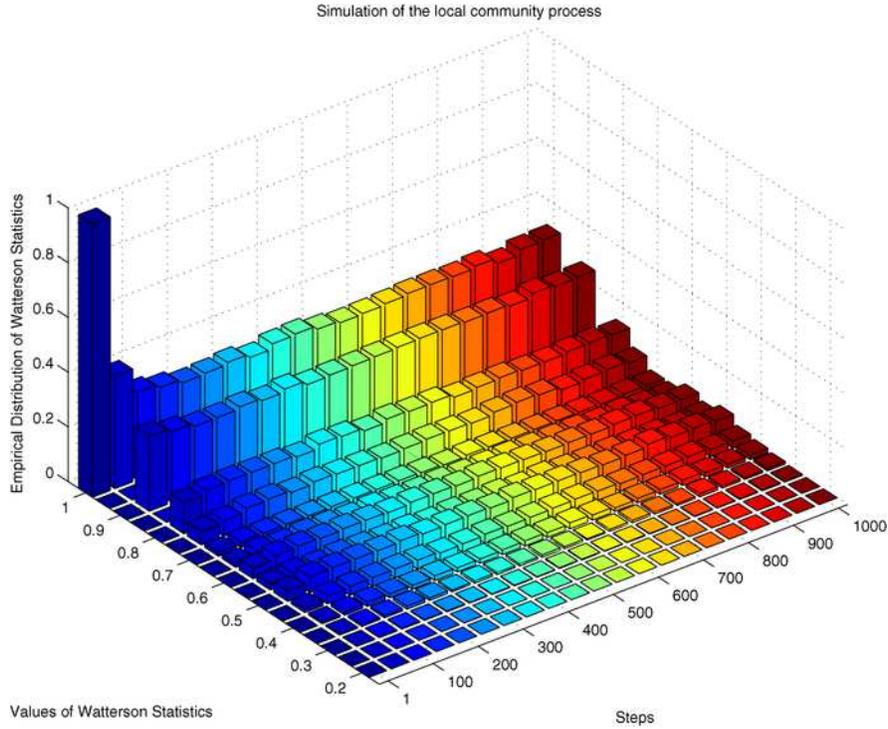

Fig. 1. *Empirical distributions of the Watterson statistics over time from 5000 simulations of the local community process. The parameters are $N = 20, d = 5, m = 0.05$ and $p_i = 1/5$. Starting state is $N\mathbf{e}_i$.*

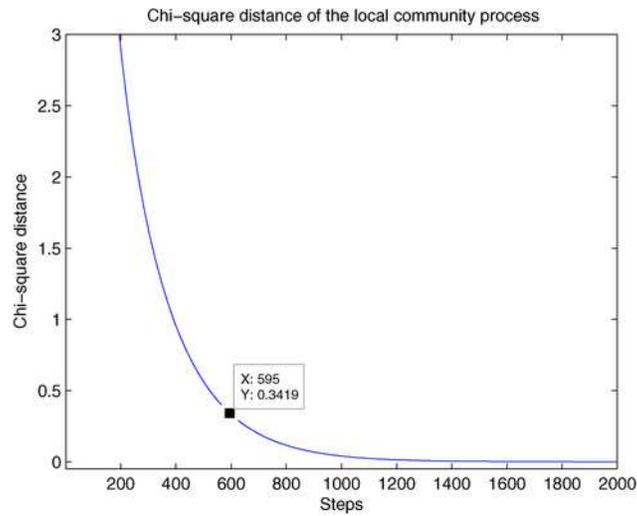

Fig. 2. *Chi-square distance of the local community process over time. The parameters are $N = 20$, $d = 5$, $m = 0.05$ and $p_i = 1/5$. Starting state is $N\mathbf{e}_i$.*

The paper is organized as follows. Section 2 provides the required background on the rates of convergence for Markov chains and some multivariate orthogonal polynomials. Section 3 gives simple criteria to verify that a reversible Markov kernel has orthogonal polynomials as eigenfunctions. Section 4 contains a wide spectrum of applications. In every example we analyze the rates of convergence of the Markov chain starting at a natural initial state. The frameworks in Sections 4.1 and 4.2 allow treatment of seemingly different Markov chains in a unified way. Specific examples, for example, the multivariate Moran model, the Dirichlet-multinomial Gibbs sampler and Bernoulli–Laplace processes extend previous results [6, 11, 15, 16] for the univariate case. Section 4.3 contains analysis of a class of Ehrenfest urn models which generalize [35]. In Section 4.4, we analyze the multivariate normal autoregressive process which arises in the multigrid Monte Carlo method [22] and general overrelaxation MCMC algorithms [3, 4, 40].

We realize that the Markov chains being actually used at the forefront of today's research are quite complicated and it is still a serious research effort to provide useful analysis for the rates of convergence of such Markov chains. Still, our examples are standard and easy to understand models and it is nice to have a sharp analysis of the rates of convergence of these chains.

## 2. Background.

2.1. *Convergence rates of Markov chains.* Let $(\mathcal{X}, \mathcal{F})$ be a measurable space equipped with a $\sigma$-finite measure $\mu$. Suppose we are given a Markov chain on state space $\mathcal{X}$ described by its transition density $K(x, x')$ with respect to $\mu(dx')$. Suppose further that the chain has stationary measure $m(dx) = m(x)\mu(dx)$. Let $K^l(x, \cdot)$ denote the density of the chain started at state $x$ after $l$ steps. The chi-square distance between $K^l(x, \cdot)$ and the stationary measure $m$ is defined by

$$\chi_x^2(l) = \int_{\mathcal{X}} \frac{[K^l(x, x') - m(x')]^2}{m(x')} \mu(dx').$$

Intuitively the chi-square distance penalizes more the discrepancies at points which have smaller probability mass (density) at stationarity. The commonly used total variation distance is defined by

$$\|K_x^l - m\|_{\mathrm{TV}} = \frac{1}{2} \int_{\mathcal{X}} |K^l(x, x') - m(x')| \mu(dx').$$

While total variation distance is always bounded in interval $[0, 1]$, the chi-square distance assumes values in $[0, \infty]$. By Cauchy–Schwarz inequality, the chi-square distance gives an upper bound for the total variation distance

$$4\|K_x^l - m\|_{\mathrm{TV}}^2 \leq \chi_x^2(l).$$

Let $l^2(m) := \{f : \mathcal{X} \to \mathbb{R} : \int_{\mathcal{X}} f^2(x) m(x) \mu(dx) < \infty\}$ denote the Hilbert space equipped with inner product

$$\langle f, g \rangle_{l^2(m)} = \mathbf{E}_m[f(X)g(X)] = \int_{\mathcal{X}} f(x)g(x)m(x)\mu(dx).$$

The Markov chain $K$ operates on $l^2(m)$ by

$$Kf(x) = \int_{\mathcal{X}} K(x,y)f(y)\mu(dy).$$

$K$ is called reversible when $m(x)K(x,x') = m(x')K(x',x)$ for all $x, x' \in \mathcal{X}$, or equivalently, when the operator $K$ is self-adjoint:

$$\langle Kf, g \rangle_{l^2(m)} = \langle f, Kg \rangle_{l^2(m)}.$$

Suppose that $l^2(m)$ admits an orthonormal basis of eigenfunctions $\{\phi_n\}_{n \geq 0}$ with $\phi_0 \equiv 1$ such that

$$K\phi_n(x) = \beta_n \phi_n(x), \qquad n \geq 0,$$

where the eigenvalues $\{\beta_n\}_{n \geq 0}$ satisfy $\beta_0 = 1$, $|\beta_n| \leq 1$, and $\sum_{n=1}^{\infty} \beta_n^2 < \infty$. Then, $K$ is a Hilbert–Schmidt operator and

$$K^l(x, x') = \sum_{n=0}^{\infty} \beta_n^l \phi_n(x) \phi_n(x') m(x') \qquad [\text{convergence in } l^2(m \times m)].$$

Also,

(2.1) $$\chi_x^2(l) = \sum_{n=1}^{\infty} \beta_n^{2l} \phi_n^2(x).$$

The central identity in our analysis throughout the paper is (2.1). The challenge is to work with the eigenvalues $\{\beta_n\}_{n \geq 0}$ and eigenfunctions $\{\phi_n\}_{n \geq 0}$ and manipulate the right-hand side in (2.1) to obtain sharp rates of convergence.

2.2. *Multivariate orthogonal polynomials.* Before introducing multivariate orthogonal polynomials we first set up the definitions of some standard multivariate distributions. To ease the notational burden, we will use boldface letters to indicate vectors. For a vector $\mathbf{x} = (x_1, \ldots, x_d) \in \mathbb{R}^d$,

$$|\mathbf{x}| = \sum_{i=1}^{d} x_i, \qquad |\mathbf{x}_i| = \sum_{j=1}^{i} x_j, \qquad |\mathbf{x}^i| = \sum_{j=i}^{d} x_j.$$

The multinomial coefficient is denoted by

$$\binom{|\mathbf{x}|}{\mathbf{x}} = \frac{|\mathbf{x}|!}{x_1! \cdots x_d!}.$$

Increasing and decreasing factorials are denoted by

$$a_{(k)} = a(a+1)\cdots(a+k-1), \qquad a_{[k]} = a(a-1)\cdots(a-k+1).$$

$a_{(0)} = a_{[0]} = 1$ by convention. For vectors $\mathbf{x} = (x_1, \ldots, x_d)$ and $\mathbf{n} = (n_1, \ldots, n_d)$,

$$\mathbf{x}^{\mathbf{n}} = \prod_{i=1}^{d} x_i^{n_i}, \qquad \mathbf{x}_{(\mathbf{n})} = \prod_{i=1}^{d} (x_i)_{(n_i)}, \qquad \mathbf{x}_{[\mathbf{n}]} = \prod_{i=1}^{d} (x_i)_{[n_i]}.$$

Three spaces that play important roles later are listed below:

$$\Delta^d = \{\mathbf{p} = (p_1, \ldots, p_d) \in [0,1]^d : |\mathbf{p}| = 1\},$$
$$\mathcal{X}_N^d = \{\mathbf{x} = (x_1, \ldots, x_d) \in \mathbb{N}_0^d : |\mathbf{x}| = N\},$$
$$\mathcal{X}_{N,\boldsymbol{\ell}}^d = \{\mathbf{x} = (x_1, \ldots, x_d) \in \mathbb{N}_0^d : |\mathbf{x}| = N, 0 \leq x_i \leq l_i\}.$$

*Dirichlet distribution* with parameters $\boldsymbol{\alpha} = (\alpha_1, \ldots, \alpha_d)$, $\alpha_i > 0$, has support $\Delta^d$ and density

$$(2.2) \qquad \mathcal{D}(\mathbf{p}|\boldsymbol{\alpha}) = \frac{\Gamma(|\boldsymbol{\alpha}|)}{\Gamma(\alpha_1)\cdots\Gamma(\alpha_d)} \prod_{i=1}^{d} p_i^{\alpha_i - 1}, \qquad \mathbf{p} \in \Delta^d.$$

In Bayesian statistics, Dirichlet distribution is the conjugate prior for the multinomial distribution.

*Multinomial distribution*, with parameters $N > 0$ and $\mathbf{p} \in \Delta^d$, has support $\mathcal{X}_N^d$ and probability mass function

$$(2.3) \qquad \mathcal{M}(\mathbf{x}|N, \mathbf{p}) = \binom{N}{\mathbf{x}} \prod_{i=1}^{d} p_i^{x_i}, \qquad \mathbf{x} \in \mathcal{X}_N^d.$$

*Dirichlet-multinomial distribution*, with parameters $N > 0$ and $\boldsymbol{\alpha} = (\alpha_1, \ldots, \alpha_d)$, $\alpha_i > 0$, is the Dirichlet $\mathcal{D}(\mathbf{p}|\boldsymbol{\alpha})$ mixture of multinomial $\mathcal{M}(\cdot|N, \mathbf{p})$ and has probability mass function

$$(2.4) \quad \mathcal{DM}(\mathbf{x}|N, \boldsymbol{\alpha}) = \binom{N}{\mathbf{x}} \frac{\prod_{i=1}^{d}(\alpha_i)_{(x_i)}}{|\boldsymbol{\alpha}|_{(N)}} = \frac{\prod_{i=1}^{d}\binom{x_i + \alpha_i - 1}{x_i}}{\binom{N + |\boldsymbol{\alpha}| - 1}{N}}, \qquad \mathbf{x} \in \mathcal{X}_N^d.$$

The same distribution is called *multivariate negative hypergeometric distribution* in [32], page 179. When $\boldsymbol{\alpha} = (1, \ldots, 1)$, it is the well-known Bose–Einstein distribution. Note as $\boldsymbol{\alpha}/|\boldsymbol{\alpha}| \to \mathbf{p} = (p_1, \ldots, p_d) \in \Delta^d$, $\mathcal{DM}(\mathbf{x}|N, \boldsymbol{\alpha}) \to \mathcal{M}(\mathbf{x}|N, \mathbf{p})$.

Let $\boldsymbol{\ell} = (l_1, \ldots, l_d)$ be a vector of positive integers and $0 < N < |\boldsymbol{\ell}|$. Replacing $\alpha_i$ by $-l_i$ for $1 \leq i \leq d$ in (2.4), we obtain the *hypergeometric distribution* with parameter $N$ and $\boldsymbol{\ell}$:

$$(2.5) \qquad \mathcal{H}(\mathbf{x}|N, \boldsymbol{\ell}) = \binom{N}{\mathbf{x}} \frac{\prod_{i=1}^{d}(l_i)_{[x_i]}}{|\boldsymbol{\ell}|_{[N]}} = \frac{\prod_{i=1}^{d}\binom{l_i}{x_i}}{\binom{|\boldsymbol{\ell}|}{N}}, \qquad \mathbf{x} \in \mathcal{X}_{N,\boldsymbol{\ell}}^d.$$

Classically the hypergeometric distribution occurs from sampling $N$ balls without replacement from a pool of $|\boldsymbol{\ell}|$ balls with composition $\boldsymbol{\ell}$.

Multinomial, Dirichlet-multinomial and hypergeometric distributions have alternative interpretations in a unified framework called the Pólya (or Pólya–Eggenberger) urn scheme. Chapter 40 of [32] is dedicated to the properties and history of the Pólya urn model. An urn initially contains $l_1$ balls of color 1, $l_2$ balls of color $2, \ldots, l_d$ balls of color $d$. In a typical Pólya urn scheme, a ball is randomly drawn. The color of the ball is noted and the ball is returned to the urn along with $c$ additional balls of the same color. This experiment is repeated $N$ times and the distribution of the composition of the observed $N$ balls is called a *Pólya–Eggenberger distribution*. When $c = 0$, this becomes *sampling with replacement* and the distribution is multinomial with parameters $N$ and $\mathbf{p} = \boldsymbol{\ell}/|\boldsymbol{\ell}|$. When $c > 0$, the distribution is Dirichlet-multinomial with parameters $N$ and $\boldsymbol{\alpha} = \boldsymbol{\ell}/c$. When $c = -1$, this becomes *sampling without replacement* and the distribution is hypergeometric with parameters $N$ and $\boldsymbol{\ell}$. In recent literature, Pólya type processes allow much more general replacement schemes. In Section 4.1, the Moran model in population genetics, its variants in community ecology, and the Dirichlet-multinomial Gibbs sampler are put in a unified framework called the sequential Pólya urn models.

The multivariate normal distribution is used when we study the multivariate normal autoregressive process in Section 4.4. The density of a *multivariate normal distribution* with mean vector $\boldsymbol{\mu} \in \mathbb{R}^d$ and covariance matrix $\Sigma$ at a column vector $\mathbf{x}$ is given by

$$\mathcal{N}(\mathbf{x}|\boldsymbol{\mu}, \Sigma) = \frac{1}{\sqrt{(2\pi)^d |\Sigma|}} e^{-(\mathbf{x}-\boldsymbol{\mu})^T \Sigma^{-1} (\mathbf{x}-\boldsymbol{\mu})/2}, \qquad \mathbf{x} \in \mathbb{R}^d.$$

Next we review the relevant developments of multivariate orthogonal polynomials for these classical multivariate distributions. These polynomials occur as the eigenfunctions for the Markov chains analyzed in this paper.

2.2.1. *Review of some explicit multivariate orthogonal polynomials.* Iliev and Xu [28] explicitly construct systems of orthogonal polynomials on various multivariate discrete weight functions which correspond to classical multivariate probability distributions. Most pertinent to us are the multivariate Hahn polynomials. We present their construction below. The parameters in [28] are shifted by one from ours.

*Multivariate Hahn polynomials.* For any $(n_1, \ldots, n_d) \in \mathcal{X}_N^d$, we let $\mathbf{n} = (n_1, \ldots, n_{d-1})$ be the index vector and define multivariate polynomials on $\mathcal{X}_N^d$ by

$$Q_{\mathbf{n}}(\mathbf{x}; N, \boldsymbol{\alpha})$$

$$= \frac{(-1)^{|\mathbf{n}|}}{(N)_{[|\mathbf{n}|]}} \prod_{j=1}^{d-1} (-N + |\mathbf{x}_{j-1}| + |\mathbf{n}^{j+1}|)_{(n_j)}$$

(2.6)
$$\times Q_{n_j}(x_j; N - |\mathbf{x}_{j-1}| - |\mathbf{n}^{j+1}|,$$
$$\alpha_j, |\boldsymbol{\alpha}^{j+1}| + 2|\mathbf{n}^{j+1}|),$$

where

$$Q_n(x; N, \alpha, \beta) = {}_3F_2\left(\begin{matrix}-n, n+\alpha+\beta-1, -x \\ \alpha, -N\end{matrix} \,\middle|\, 1\right)$$
$$= \sum_{j=0}^n \frac{(-n)_{(j)}(n+\alpha+\beta-1)_{(j)}(-x)_{(j)}}{(\alpha)_{(j)}(-N)_{(j)} j!}$$

is the classical univariate Hahn polynomial. The univariate Hahn polynomials satisfy the orthogonality relation

(2.7)
$$\mathbf{E}_{\mathcal{DM}(\cdot|N,\alpha,\beta)}[Q_n(X; N, \alpha, \beta) Q_m(X; N, \alpha, \beta)]$$
$$= \frac{(N + \alpha + \beta)_{(n)}(\beta)_{(n)}}{\binom{N}{n}(2n + \alpha + \beta - 1)(\alpha + \beta)_{(n-1)}(\alpha)_{(n)}} \delta_{mn}.$$

A survey of the univariate Hahn polynomials is in [29], Section 6.2. The following proposition is essentially Theorem 5.4 in [28].

PROPOSITION 2.1. *The system (2.6) satisfies the orthogonality relation*

(2.8)
$$\mathbf{E}_{\mathcal{DM}(\cdot|N,\boldsymbol{\alpha})}[Q_{\mathbf{n}}(\mathbf{X}; N, \boldsymbol{\alpha}) Q_{\mathbf{n}'}(\mathbf{X}; N, \boldsymbol{\alpha})]$$
$$= \sum_{\mathbf{x} \in \mathcal{X}_N^d} Q_{\mathbf{n}}(\mathbf{x}; N, \boldsymbol{\alpha}) Q_{\mathbf{n}'}(\mathbf{x}; N, \boldsymbol{\alpha}) \mathcal{DM}(\mathbf{x}|N, \boldsymbol{\alpha})$$
$$= d_{\mathbf{n}}^2 \delta_{\mathbf{n},\mathbf{n}'},$$

*where*

(2.9)
$$d_{\mathbf{n}}^2 = \frac{(|\boldsymbol{\alpha}| + N)_{(|\mathbf{n}|)}}{(N)_{[|\mathbf{n}|]} |\boldsymbol{\alpha}|_{(2|\mathbf{n}|)}}$$
$$\times \prod_{j=1}^{d-1} \frac{(|\boldsymbol{\alpha}^j| + |\mathbf{n}^j| + |\mathbf{n}^{j+1}| - 1)_{(n_j)}(|\boldsymbol{\alpha}^{j+1}| + 2|\mathbf{n}^{j+1}|)_{(n_j)} n_j!}{(\alpha_j)_{(n_j)}}.$$

In probabilistic language, the construction of (2.6) uses the stick-breaking property of the Dirichlet-multinomial distribution. Essentially the same result was implied in an earlier work by Karlin and McGregor [36].

REMARK 2.2. Here are a few properties of the orthogonal system (2.6):

1. $Q_\mathbf{n}$ is a polynomial in $x_1, \ldots, x_{d-1}$ of degree $|\mathbf{n}|$.
2. $Q_\mathbf{n}(N\mathbf{e}_d; N, \boldsymbol{\alpha}) = 1$.
3. When $d = 2$, we recover the univariate Hahn polynomials.

The orthogonality relation (2.8) does not depend on the positivity of $\alpha_i, 1 \leq i \leq d$. Therefore switching parameter range gives the orthogonal polynomials on the multivariate hypergeometric distribution. Note that the sample space (and thus the index set) is $\mathcal{X}^d_{N,\boldsymbol{\ell}}$, a subset of $\mathcal{X}^d_N$.

PROPOSITION 2.3.　*The system (2.6), with $\alpha_i = -l_i, 1 \leq i \leq d$, is orthogonal with respect to the hypergeometric distribution $\mathcal{H}(\cdot|N, \boldsymbol{\ell})$ (2.5).*

*Multivariate Krawtchouk polynomials.*　Recall the well-known limiting relation

$$Q_n(x; N, pt, (1-p)t) = \sum_{j=0}^{n} \frac{(-n)_{(j)}(n + t - 1)_{(j)}(-x)_j}{(pt)_{(j)}(-N)_{(j)} j!}$$

$$\to \sum_{j=0}^{n} \frac{(-n)_{(j)}(-x)_{(j)}}{(-N)_{(j)} j! p^j}$$

$$= {}_2F_1\left(\begin{matrix}-n, -x \\ -N\end{matrix} \,\bigg|\, \frac{1}{p}\right) = K_n(x; N, p)$$

as $t \to \infty$, where $K_n$'s are the univariate Krawtchouk polynomials orthogonal for the binomial distribution with parameters $N$ and $p$. Properties of the univariate Krawtchouk polynomials are documented in [29], page 183. Taking limits $\boldsymbol{\alpha}/|\boldsymbol{\alpha}| \to \mathbf{p} = (p_1, \ldots, p_d) \in \Delta^d$ in the orthogonality relation of the multivariate Hahn polynomials (2.8), we obtain the orthogonal polynomials (multivariate Krawtchouk) for the multinomial distribution $\mathcal{M}(\cdot|N, \mathbf{p})$.

PROPOSITION 2.4.　*The system*

(2.10)
$$K_\mathbf{n}(\mathbf{x}; N, \mathbf{p}) = \frac{(-1)^{|\mathbf{n}|}}{(N)_{[|\mathbf{n}|]}} \prod_{j=1}^{d-1} (-N + |\mathbf{x}_{j-1}| + |\mathbf{n}^{j+1}|)_{(n_j)}$$

$$\times K_{n_j}\left(x_j; N - |\mathbf{x}_{j-1}| - |\mathbf{n}^{j+1}|, \frac{p_j}{|\mathbf{p}^j|}\right)$$

*satisfies the orthogonality relation*

$$\mathbf{E}_{\mathcal{M}(\cdot|N,\mathbf{p})}[K_\mathbf{n}(\mathbf{X}; N, \mathbf{p}) K_{\mathbf{n}'}(\mathbf{X}; N, \mathbf{p})]$$

$$= \sum_{\mathbf{x} \in \mathcal{X}^d_N} K_\mathbf{n}(\mathbf{x}; N, \mathbf{p}) K_{\mathbf{n}'}(\mathbf{x}; N, \mathbf{p}) \mathcal{M}(\mathbf{x}|N, \mathbf{p})$$

$$= d^2_\mathbf{n} \delta_{\mathbf{n},\mathbf{n}'},$$

*where*

$$d_{\mathbf{n}}^2 = \frac{1}{(N)_{[|\mathbf{n}|]}} \prod_{j=1}^{d-1} \frac{(|\mathbf{p}^j|)^{n_j}(|\mathbf{p}^{j+1}|)^{n_j} n_j!}{p_j^{n_j}}. \tag{2.11}$$

The systems of orthogonal polynomials for a fixed multinomial distribution are not unique. A more general construction of the multivariate Krawtchouk polynomials was given by Griffiths [23] and recently surveyed in [43]. The system defined by (2.10) is a special case in this general framework.

*Multivariate Jacobi polynomials.* Next we record the multivariate orthogonal polynomials for the Dirichlet distribution. Recall that, for the univariate Hahn polynomials,

$$\begin{aligned}
Q_n(x; N, \alpha, \beta) &= {}_3F_2\left(\begin{matrix}-n, n+\alpha+\beta-1, -x \\ \alpha, -N\end{matrix} \,\middle|\, 1\right) \\
&= \sum_{j=0}^n \frac{(-n)_{(j)}(n+\alpha+\beta-1)_{(j)}(-x)_{(j)}}{(\alpha)_{(j)}(-N)_{(j)} j!} \\
&\to \sum_{j=0}^n \frac{(-n)_{(j)}(n+\alpha+\beta-1)_{(j)}}{(\alpha)_{(j)}} \frac{z^j}{j!} \\
&= {}_2F_1\left(\begin{matrix}-n, n+\alpha+\beta-1 \\ \alpha\end{matrix} \,\middle|\, z\right) = J_n(z; \alpha, \beta)
\end{aligned}$$

as $x/N \to z$. $\{J_n\}_{0 \le n < \infty}$ are the shifted Jacobi polynomials which are orthogonal for the beta distribution with parameters $\alpha$ and $\beta$. Taking limits $\mathbf{x}/N \to \mathbf{z} = (z_1, \ldots, z_d) \in \Delta^d$ in the multivariate Hahn polynomials (2.6), we obtain a system of multivariate polynomials on $\Delta^d$:

$$J_{\mathbf{n}}(\mathbf{z}; \boldsymbol{\alpha}) = \prod_{j=1}^{d-1} |\mathbf{z}^j|^{n_j} J_{n_j}\left(\frac{z_j}{|\mathbf{z}^j|}; \alpha_j, |\boldsymbol{\alpha}^{j+1}| + 2|\mathbf{n}^{j+1}|\right), \tag{2.12}$$

$$\mathbf{z} \in \Delta^d, \ \mathbf{n} \in \mathbb{N}_0^{d-1}.$$

PROPOSITION 2.5. *The system (2.12) satisfies the orthogonality relation*

$$\begin{aligned}
\mathbf{E}_{\mathcal{D}(\cdot|\boldsymbol{\alpha})}&[J_{\mathbf{n}}(\mathbf{Z}; \boldsymbol{\alpha}) J_{\mathbf{n}'}(\mathbf{Z}; \boldsymbol{\alpha})] \\
&= \int_{\Delta^d} J_{\mathbf{n}}(\mathbf{z}; \boldsymbol{\alpha}) J_{\mathbf{n}'}(\mathbf{z}; \boldsymbol{\alpha}) \mathcal{D}(\mathbf{z}|\boldsymbol{\alpha}) \, d\mathbf{z} \\
&= d_{\mathbf{n}}^2 \delta_{\mathbf{n}, \mathbf{n}'},
\end{aligned}$$

*where*

$$d_{\mathbf{n}}^2 = \frac{1}{|\boldsymbol{\alpha}|_{(2|\mathbf{n}|)}}$$

(2.13)
$$\times \prod_{j=1}^{d-1} \frac{(|\boldsymbol{\alpha}^j| + |\mathbf{n}^j| + |\mathbf{n}^{j+1}| - 1)_{(n_j)}(|\boldsymbol{\alpha}^{j+1}| + 2|\mathbf{n}^{j+1}|)_{(n_j)} n_j!}{(\alpha_j)_{(n_j)}}.$$

A stick-breaking construction of (2.12) was known earlier (see, e.g., [17, 37]).

*Hermite polynomials.* We require univariate Hermite polynomials to carry out the analysis of the multivariate normal autoregressive process. The eigenfunctions of this process (after appropriate variable transformations) turn out to be products of appropriate univariate Hermite polynomials. The univariate Hermite polynomials are defined by

(2.14) $$H_n(x) = n! \sum_{k=0}^{[n/2]} \frac{(-1)^k (2x)^{n-2k}}{k!(n-2k)!}, \qquad n \geq 0.$$

They satisfy the orthogonality relation

$$\frac{1}{\sqrt{\pi}} \int_{-\infty}^{\infty} H_m(y) H_n(y) e^{-y^2} \, dy = 2^n n! \delta_{mn}.$$

An important multilinear generating function formula ([29], Example 4.7.3), gives

(2.15) $$\sum_{n=0}^{\infty} \frac{H_n^2(x)}{2^n n!} t^n = \frac{1}{\sqrt{1-t^2}} e^{2x^2 t/(1+t)}.$$

2.2.2. *Griffiths' construction of kernel polynomials.* In the applications in Section 4, we need to manipulate certain sums of products of the relevant multivariate orthogonal polynomials. For a fixed multivariate distribution $m$, the *kernel polynomials* are defined as

$$h_n(\mathbf{x}, \mathbf{y}) = \sum_{|\mathbf{n}|=n} Q_{\mathbf{n}}^0(\mathbf{x}) Q_{\mathbf{n}}^0(\mathbf{y})$$

for *any* complete system of orthonormal polynomials $\{Q_{\mathbf{n}}^0\}$ in $l^2(m)$. The kernel polynomials are invariant under the choice of orthogonal polynomial system. Fortunately, these sums can be carried out in closed form for the systems of polynomials that we consider. In this section, we review the work by Griffiths, who explicitly constructed the kernel polynomials for the Dirichlet-multinomial [25], Dirichlet [24, 25] and the multinomial [26] distributions.

PROPOSITION 2.6.  *For $\mathbf{x}, \mathbf{y} \in \mathcal{X}_N^d$ and $0 \leq n \leq N$, the kernel polynomials for the Dirichlet-multinomial distribution $\mathcal{DM}(\cdot|N, \boldsymbol{\alpha})$ (2.4) are*

$$h_n(\mathbf{x}, \mathbf{y}) = (|\boldsymbol{\alpha}| + 2n - 1)\frac{(|\boldsymbol{\alpha}| + N)_{(n)}}{N_{[n]}}$$

(2.16)
$$\times \sum_{m=0}^{n} (-1)^{n-m} \frac{(|\boldsymbol{\alpha}| + m)_{(n-1)}}{m!(n-m)!} \xi_m(\mathbf{x}, \mathbf{y}),$$

*where*

$$(2.17) \quad \xi_m(\mathbf{x}, \mathbf{y}) = \sum_{|\boldsymbol{\ell}|=m} \binom{m}{\boldsymbol{\ell}} \frac{|\boldsymbol{\alpha}|_{(m)}}{\prod_1^d \alpha_{i(l_i)}} \frac{\prod_1^d (\alpha_i + x_i)_{(l_i)} (\alpha_i + y_i)_{(l_i)}}{(|\boldsymbol{\alpha}| + N)_{(m)} (|\boldsymbol{\alpha}| + N)_{(m)}}.$$

For example, the first two kernel polynomials are

$$h_0(\mathbf{x}, \mathbf{y}) \equiv 1;$$

$$h_1(\mathbf{x}, \mathbf{y}) = \frac{(|\boldsymbol{\alpha}| + 1)(|\boldsymbol{\alpha}| + N)}{N} \left( \frac{|\boldsymbol{\alpha}|}{(|\boldsymbol{\alpha}| + N)^2} \sum_{i=1}^d \frac{(\alpha_i + x_i)(\alpha_i + y_i)}{\alpha_i} - 1 \right).$$

When calculating convergence rates of Markov chain starting from state $N\mathbf{e}_i$, we need to evaluate the quantity

(2.18)
$$h_n(N\mathbf{e}_i, N\mathbf{e}_i) = \frac{N_{[n]}}{(N + |\boldsymbol{\alpha}|)_{(n)}} \frac{(|\boldsymbol{\alpha}| + 2n - 1)(|\boldsymbol{\alpha}|)_{(n-1)}(|\boldsymbol{\alpha}| - \alpha_i)_{(n)}}{n!(\alpha_i)_{(n)}}$$

$$= \binom{N}{n} \frac{(|\boldsymbol{\alpha}| + 2n - 1)(|\boldsymbol{\alpha}|)_{(n-1)}(|\boldsymbol{\alpha}| - \alpha_i)_{(n)}}{(|\boldsymbol{\alpha}| + N)_{(n)}(\alpha_i)_{(n)}}.$$

The first equality follows from (2.9), (2.13) and (2.21) below.

Replacing $\alpha_i$ by $-l_i$, $1 \leq i \leq d$, in (2.16) and (2.17), we obtain the kernel polynomials for the hypergeometric distribution $\mathcal{H}(\cdot|N, \boldsymbol{\ell})$ (2.5). For example, when $l_i \geq N$ and $\mathbf{x} = \mathbf{y} = N\mathbf{e}_i$,

$$h_n(N\mathbf{e}_i, N\mathbf{e}_i) = \binom{N}{n} \frac{(|\boldsymbol{\ell}| - 2n + 1)|\boldsymbol{\ell}|_{[n-1]}(|\boldsymbol{\ell}| - l_i)_{[n]}}{(|\boldsymbol{\ell}| - N)_{[n]}(l_i)_{[n]}}.$$

Writing $x_i = Nw_i$, $y_i = Nz_i$ and taking limit $N \to \infty$ in (2.17), we recover the kernel polynomials for the Dirichlet distribution.

PROPOSITION 2.7.  *For $\mathbf{w}, \mathbf{z} \in \Delta^d$ and $0 \leq n < \infty$, the kernel polynomials for the Dirichlet distribution $\mathcal{D}(\cdot|\boldsymbol{\alpha})$ (2.2) are*

$$(2.19) \ \ h_n(\mathbf{w}, \mathbf{z}) = (|\boldsymbol{\alpha}| + 2n - 1) \sum_{m=0}^{n} (-1)^{n-m} \frac{(|\boldsymbol{\alpha}| + m)_{(n-1)}}{m!(n-m)!} \xi_m(\mathbf{w}, \mathbf{z}),$$

where

$$\xi_m(\mathbf{w}, \mathbf{z}) = \sum_{|\boldsymbol{\ell}|=m} \binom{m}{\boldsymbol{\ell}} \frac{|\boldsymbol{\alpha}|_{(m)}}{\prod_{i=1}^d \alpha_{i(l_i)}} \prod_{i=1}^d (w_i z_i)^{l_i}. \tag{2.20}$$

The kernel polynomials for the Dirichlet distribution were derived in [24] where the transition density of the multi-allele Wright–Fisher diffusion process was expanded in terms of (2.19). When $\mathbf{w} = \mathbf{z} = \mathbf{e}_i$,

$$\begin{aligned}
h_n(\mathbf{e}_i, \mathbf{e}_i) &= \frac{(|\boldsymbol{\alpha}|+2n-1)}{n!} \sum_{m=0}^n (-1)^{n-m} \binom{n}{m} \frac{|\boldsymbol{\alpha}|_{(m+n-1)}}{(\alpha_i)_{(m)}} \\
&= \frac{(|\boldsymbol{\alpha}|+2n-1)}{n!} \frac{|\boldsymbol{\alpha}|_{(n)}}{|\boldsymbol{\alpha}|+n-1} (-1)^n \sum_{m=0}^n \frac{(-n)_{(m)}(|\boldsymbol{\alpha}|+n-1)_{(m)}}{(\alpha_i)_{(m)} m!} \\
&= \frac{(|\boldsymbol{\alpha}|+2n-1)(|\boldsymbol{\alpha}|)_{(n-1)}(|\boldsymbol{\alpha}|-\alpha_i)_{(n)}}{n!(\alpha_i)_{(n)}},
\end{aligned} \tag{2.21}$$

where the last equality uses the Chu–Vandermonde summation formula.

The kernel polynomials for the multinomial distribution were given in the work by Griffiths [26].

PROPOSITION 2.8. *For $\mathbf{x}, \mathbf{y} \in \mathcal{X}_N^d$ and $0 \le n \le N$, the kernel polynomials for the multinomial distribution $\mathcal{M}(\cdot|N, \mathbf{p})$ (2.3) are*

$$h_n(\mathbf{x}, \mathbf{y}) = \sum_{m=0}^n \binom{N}{m} \binom{N-m}{n-m} (-1)^{n-m} \xi_m, \tag{2.22}$$

where

$$\xi_m = \sum_{|\boldsymbol{\ell}|=m} \binom{m}{\boldsymbol{\ell}} \frac{\prod_{i=1}^d (x_i)_{[l_i]} (y_i)_{[l_i]} p_i^{-l_i}}{N_{[m]} N_{[m]}}.$$

When $\mathbf{x} = \mathbf{y} = N\mathbf{e}_i$,

$$h_n(N\mathbf{e}_i, N\mathbf{e}_i) = \binom{N}{n} \left(\frac{1-p_i}{p_i}\right)^n. \tag{2.23}$$

**3. Markov chains with orthogonal polynomial eigenfunctions.** There is a large class of Markov chains with polynomial eigenfunctions, for example, birth–death processes [41], Cannings exchangeable model [5] in population genetics, certain two-component Gibbs samplers [11] and of course all Markov chains considered in this paper. When the Markov kernel is reversible, often the orthogonal polynomials for the stationary distribution come up as eigenfunctions. We present simple tools for identifying cases when orthogonal polynomials are the eigenfunctions of a reversible Markov kernel, first in the univariate case and then the multivariate case.

LEMMA 3.1. *Suppose $\pi$ is a univariate distribution and $l^2(\pi)$ admits an orthogonal basis of polynomials $\{q_n(x)\}_{0 \leq n < c}$ where $c = \#\mathrm{supp}(\pi)$, that is, $q_n(x)$ is a polynomial in $x$ of exact degree $n$ and*

$$\langle q_n, q_m \rangle_{l^2(\pi)} = \mathbf{E}_\pi[q_n(X)q_m(X)] = d_n^2 \delta_{nm}.$$

*If $K$ is a Markov kernel reversible on $\pi$ and*

(3.1) $\quad \mathbf{E}_{K(x,\cdot)}[X^n] = \beta_n x^n + \text{terms in } x \text{ of degree} < n, \qquad 0 \leq n < c,$

*then:*

1. *$K$ has eigenvalue $\beta_n$ with eigenfunction $q_n$, that is, $Kq_n = \beta_n q_n$ for $0 \leq n < c$.*
2. *The chi-square distance between $K^l(x, \cdot)$ and $\pi$ is*

$$\chi_x^2(l) = \sum_{n \geq 1} \beta_n^{2l} q_n^2(x) d_n^{-2}.$$

PROOF. $q_0$ is a constant function and trivially $Kq_0 = q_0$ with eigenvalue $\beta_0 = 1$. By hypothesis (3.1), $Kq_n = \beta_n q_n + \sum_{m<n} a_m q_m$. But the coefficients $a_m = \langle Kq_n, q_m \rangle_{l^2(\pi)} = \langle q_n, Kq_m \rangle_{l^2(\pi)} = 0$ since $Kq_m$ is a polynomial of degree $m$ and can be expanded in terms of the basis polynomials of degree $\leq m < n$ which are all orthogonal to $q_n$. Therefore $Kq_n = \beta_n q_n$. The expression for the chi-square distance follows from (2.1). □

Under mild assumptions, this result can be generalized into a multivariate case. We recall that the notation $|\mathbf{n}|$ denotes the sum of coordinates of a vector $\mathbf{n}$.

LEMMA 3.2. *Suppose $\pi$ is a multivariate distribution and $l^2(\pi)$ admits an orthogonal basis of multivariate polynomials $\{q_\mathbf{n}(\mathbf{x})\}$ where $q_\mathbf{n}$ is a polynomial of exact degree $|\mathbf{n}|$ and*

$$\langle q_\mathbf{n}, q_\mathbf{m} \rangle_{l^2(\pi)} = \mathbf{E}_\pi[q_\mathbf{n}(X)q_\mathbf{m}(X)] = d_\mathbf{n}^2 \delta_{\mathbf{nm}}.$$

*If $K$ is a Markov kernel reversible with respect to $\pi$ and*

(3.2) $\qquad \mathbf{E}_{K(\mathbf{x},\cdot)}[\mathbf{X}^\mathbf{n}] = \beta_{|\mathbf{n}|} \mathbf{x}^\mathbf{n} + \text{terms in } \mathbf{x} \text{ of degree} < |\mathbf{n}|,$

*then:*

1. *$K$ has eigenvalue $\beta_n$ with corresponding eigenbasis $\{q_\mathbf{n}\}_{|\mathbf{n}|=n}$.*
2. *The chi-square distance between $K^l(\mathbf{x}, \cdot)$ and $\pi$ is*

$$\chi_\mathbf{x}^2(l) = \sum_{n \geq 1} \beta_n^{2l} h_n(\mathbf{x}, \mathbf{x}),$$

*where*

$$h_n(\mathbf{x}, \mathbf{y}) := \sum_{|\mathbf{n}|=n} q_\mathbf{n}(\mathbf{x}) q_\mathbf{n}(\mathbf{y}) d_\mathbf{n}^{-2}$$

*is the kernel polynomial of degree n for* $\pi$.

PROOF. $q_\mathbf{0}$ is a constant function and trivially $Kq_\mathbf{0} = q_\mathbf{0}$ with eigenvalue $\beta_0 = 1$. By hypothesis (3.2), $Kq_\mathbf{n} = \beta_{|\mathbf{n}|} q_\mathbf{n} + \sum_{|\mathbf{m}|<|\mathbf{n}|} a_\mathbf{m} q_\mathbf{m}$. But the coefficients $a_\mathbf{m} = \langle Kq_\mathbf{n}, q_\mathbf{m} \rangle_{l^2(\pi)} = \langle q_\mathbf{n}, Kq_\mathbf{m} \rangle_{l^2(\pi)} = 0$ since $Kq_\mathbf{m}$ is a polynomial of degree $|\mathbf{m}|$ and can be expanded in terms of the basis polynomials of degree $\leq |\mathbf{m}| < |\mathbf{n}|$, which are all orthogonal to $q_\mathbf{n}$. Therefore $Kq_\mathbf{n} = \beta_{|\mathbf{n}|} q_\mathbf{n}$. The expression for the chi-square distance follows from (2.1). □

REMARK 3.3. When checking conditions (3.1) or (3.2), often it is easier to calculate the factorial moments. Because of the simple relation $x_{[n]} = \sum_{k=0}^{n} s(n,k) x^k$, where $s(n,k)$ are the Stirling numbers of the first kind and especially $s(n,n) = 1$, the condition (3.2) is equivalent to

$$\mathbf{E}_{K(\mathbf{x},\cdot)}[\mathbf{X}_{[\mathbf{n}]}] = \beta_{|\mathbf{n}|} \mathbf{x}_{[\mathbf{n}]} + \text{terms in } \mathbf{x} \text{ of degree} < |\mathbf{n}|.$$

The key condition (3.2) seems restrictive but holds for surprisingly many multivariate Markov chains which possess certain intrinsic symmetry. Most examples in this paper satisfy (3.2). See [43] for a class of multinomial chains with orthogonal polynomial eigenfunctions but which in general do not satisfy (3.2).

The trick of preserving polynomials has a long history in population genetics (see, for example, [5, 19]). But most models in population genetics, for example, Wright–Fisher model, are irreversible and thus orthogonal polynomials do not come up. An exception is the Moran process which we study and generalize in Section 4.1.

## 4. Applications.

4.1. *Sequential Pólya urn models.* In this section, we study a class of multivariate Markov chains which can be described in terms of the classical Pólya urns. They are all reversible with respect to the Dirichlet-multinomial distribution and have the multivariate Hahn polynomials as eigenfunctions. Interestingly, the classical multi-allele Moran process in population genetics, local community process in community ecology, and the Dirichlet-multinomial Gibbs sampler in statistics can be treated as special cases in this unified framework. The urn description also provides a convenient way for simulation of the processes on a computer.

We first define the Pólya type urns. A newly constituted urn contains one ball of color $i$ and weight $\alpha_i$ for $1 \leq i \leq d$, that is, total of $d$ balls with total weight $|\boldsymbol{\alpha}|$. A *Pólya type draw* is defined as a random draw according to weights and the ball is returned to the urn along with one additional ball of the same color and of unit weight. Initially a batch of $N$ balls of composition $(X_{01}, \ldots, X_{0d})$, that is, $X_{0i}$ balls of color $i$ and each of unit weight, are added into a newly constituted urn. We define three classes of Pólya urn models. In the following, $s$ is a fixed integer between 0 and $N$.

DEFINITION 4.1. Pólya level models: One step of the Markov chain includes three mini-steps. First randomly mark $s$ balls among the $N$ balls (excluding the $d$ original balls) to be removed. Before removal, do $s$ Pólya type draws (including the $d$ original balls). Then the $s$ marked balls are removed. Apparently the number of balls in the urn is kept constant after each step and the $d$ original balls are always in the urn. Let $\mathbf{X}_t = (X_{t1}, \ldots, X_{td})$ be the composition of $N$ balls (excluding the $d$ original balls) after $t$ steps. $\{\mathbf{X}_t\}_{t \geq 0}$ forms a multivariate Markov chain on $\mathcal{X}_N^d$, which we call a Pólya level model.

DEFINITION 4.2. Pólya down-up models: one-step of the Markov chain includes two mini-steps. First randomly choose $s$ balls among the $N$ balls (excluding the $d$ original balls) and remove them. Then do $s$ Pólya type draws (including the $d$ original balls). The number of balls in the urn is kept constant and the $d$ original balls are always in the urn. Let $\mathbf{X}_t = (X_{t1}, \ldots, X_{td})$ be the composition of $N$ balls (excluding the $d$ original balls) after $t$ steps. $\{\mathbf{X}_t\}_{t \geq 0}$ forms a multivariate Markov chain on $\mathcal{X}_N^d$, which we call a Pólya down-up model.

DEFINITION 4.3. Pólya up-down models: one-step of the Markov chain includes two mini-steps. First do $s$ Pólya type draws (including the $d$ original balls). Then randomly choose $s$ balls among the $N+s$ balls (excluding the $d$ original balls) and remove them. The number of balls in the urn is kept constant and the $d$ original balls are always in the urn. Let $\mathbf{X}_t = (X_{t1}, \ldots, X_{td})$ be the composition of $N$ balls (excluding the $d$ original balls) after $t$ steps. $\{\mathbf{X}_t\}_{t \geq 0}$ forms a multivariate Markov chain on $\mathcal{X}_N^d$, which we call a Pólya up-down model.

All three classes of the sequential Pólya models have the same stationary distribution. The detailed balance is checked in a straightforward way.

LEMMA 4.4. *All three classes of Pólya urn models are reversible with respect to the Dirichlet-multinomial distribution $\mathcal{DM}(\cdot|N, \boldsymbol{\alpha})$ (2.4), where $\boldsymbol{\alpha} = (\alpha_1, \ldots, \alpha_d)$.*

Using the trick of preserving polynomials (Lemma 3.2), it is easy to check that all three classes of Pólya urn models take the multivariate Hahn polynomials (2.6) as eigenfunctions. Recall that $|\boldsymbol{\alpha}| = \sum_{i=1}^{d} \alpha_i$.

THEOREM 4.5. 1. *The Pólya level model has eigenvalue*

$$\begin{aligned}
\beta_n &= \sum_{k=0}^{n} \binom{n}{k} \frac{(N-s)_{[k]} s_{[n-k]}}{N_{[k]}(N+|\boldsymbol{\alpha}|)_{(n-k)}} \\
&= \sum_{k=0}^{n} \binom{n}{k} \frac{(N-s)_{[n-k]} s_{[k]}}{N_{[n-k]}(N+|\boldsymbol{\alpha}|)_{(k)}}, \qquad 0 \leq n \leq N,
\end{aligned} \tag{4.1}$$

*with multiplicity $\binom{d+n-2}{n}$ and corresponding eigenbasis $\{Q_\mathbf{n}\}_{|\mathbf{n}|=n}$, where $Q_\mathbf{n}$ are the multivariate Hahn polynomials (2.6).*

2. *The Pólya down-up model has eigenvalue*

$$\beta_n = \frac{(N-s)_{[n]}(N+|\boldsymbol{\alpha}|)_{(n)}}{N_{[n]}(N-s+|\boldsymbol{\alpha}|)_{(n)}}, \qquad 0 \leq n \leq N,$$

*with multiplicity $\binom{d+n-2}{n}$ and corresponding eigenbasis $\{Q_\mathbf{n}\}_{|\mathbf{n}|=n}$, where $Q_\mathbf{n}$ are the multivariate Hahn polynomials (2.6).*

3. *The Pólya up-down model has eigenvalue*

$$\beta_n = \frac{N_{[n]}(N+s+|\boldsymbol{\alpha}|)_{(n)}}{(N+s)_{[n]}(N+|\boldsymbol{\alpha}|)_{(n)}}, \qquad 0 \leq n \leq N,$$

*with multiplicity $\binom{d+n-2}{n}$ and corresponding eigenbasis $\{Q_\mathbf{n}\}_{|\mathbf{n}|=n}$, where $Q_\mathbf{n}$ are the multivariate Hahn polynomials (2.6).*

4. *For all three classes of Pólya models, the chi-square distance between $K^l(\mathbf{x}, \cdot)$ and the stationary distribution is*

$$\chi^2_\mathbf{x}(l) = \sum_{n=1}^{N} \beta_n^{2l} h_n(\mathbf{x}, \mathbf{x}), \tag{4.2}$$

*where $\beta_n$ are the eigenvalues for the corresponding process and $h_n$ are the kernel polynomials (2.16) for the Dirichlet-multinomial distribution.*

PROOF. For sake of space, we only provide the proof for the Pólya level models. Proofs for the other models are similar. For the Pólya level model with parameter $s$, it is observed that, given $\mathbf{X}_t = \mathbf{x}$, $\mathbf{X}_{t+1} = \mathbf{x} - \mathbf{Y} + \mathbf{Z}$, where $\mathbf{Y}$ is multivariate hypergeometric $\mathcal{H}(\cdot|s, \mathbf{x})$, $\mathbf{Z}$ is Dirichlet-multinomial $\mathcal{DM}(\cdot|s, \boldsymbol{\alpha} + \mathbf{x})$, and $\mathbf{Y}$ is independent of $\mathbf{Z}$. Joint factorial moments of multivariate hypergeometric and Dirichlet-multinomial distributions are well known. With $\mathbf{x} = (x_1, \ldots, x_{d-1})$ and $\mathbf{n} = (n_1, \ldots, n_{d-1})$,

$$\mathbf{E}_{K(\mathbf{x},\cdot)} \mathbf{X}_{[\mathbf{n}]}$$
$$= \mathbf{E}(\mathbf{x} - \mathbf{Y} + \mathbf{Z})_{[\mathbf{n}]}$$
$$= \sum_{\mathbf{0} \leq \mathbf{k} \leq \mathbf{n}} \prod_{j=1}^{d-1} \binom{n_j}{k_j} \mathbf{E}(\mathbf{x} - \mathbf{Y})_{[\mathbf{n}-\mathbf{k}]} \mathbf{E} \mathbf{Z}_{[\mathbf{k}]}$$
$$= \sum_{\mathbf{0} \leq \mathbf{k} \leq \mathbf{n}} \prod_{j=1}^{d-1} \binom{n_j}{k_j} \frac{(N-s)_{[|\mathbf{n}|-|\mathbf{k}|]}}{N_{[|\mathbf{n}|-|\mathbf{k}|]}} \mathbf{x}_{[\mathbf{n}-\mathbf{k}]} \cdot \frac{s_{[|\mathbf{k}|]}}{(N+|\boldsymbol{\alpha}|)_{(|\mathbf{k}|)}} (\mathbf{x} + \boldsymbol{\alpha})_{(\mathbf{k})}$$
$$= \left[ \sum_{k=0}^{|\mathbf{n}|} \binom{|\mathbf{n}|}{k} \frac{(N-s)_{[|\mathbf{n}|-k]} s_{[k]}}{N_{[|\mathbf{n}|-k]}(N+|\boldsymbol{\alpha}|)_{(k)}} \right] \mathbf{x}_{[\mathbf{n}]} + \text{terms in } \mathbf{x} \text{ of degree} < |\mathbf{n}|.$$

Then the claims follow from Lemma 3.2. □

REMARK 4.6. For all three classes of Pólya urn models, the eigenvalues depend on the parameter $\boldsymbol{\alpha} = (\alpha_1, \ldots, \alpha_d)$ only through its sum $|\boldsymbol{\alpha}|$. Formulas for the eigenfunctions require complete knowledge of $\boldsymbol{\alpha}$.

Complete spectral information allows for sharp results in convergence rate in chi-square distance. In case $s = 1$, an argument similar to that in [21] can also be used to obtain convergence rate in separation distance using only eigenvalues (see [44]). But in the current paper we confine ourselves to the convergence rate in chi-square distance. The following three examples are special cases of the sequential Pólya urn models.

4.1.1. *Convergence rate of the Moran process in population genetics.* In brief, the classical Moran process in population genetics models the evolution of a population of constant size by random replacement followed by mutation. Suppose there are $d$ species in a population of size $N$. At each step, one individual is chosen uniformly to die and independently another is chosen uniformly to reproduce. *They may be the same individual.* If the latter is of species $i$, the offspring has probability $m_{ij}$, $1 \leq j \leq d$, to mutate to type $j$. Let $\mathbf{X}_t = (X_{t1}, \ldots, X_{td})$ be the vector of counts of species $1, \ldots, d$ at time $t$. $\{\mathbf{X}_t\}_{t \geq 0}$ forms a Markov chain on $\mathcal{X}_N^d$. The size of the state space is $|\mathcal{X}_N^d| = \binom{N+d-1}{N}$. Let $\mathbf{x} \in \mathcal{X}_N^d$; one-step transition probabilities are

$$K(\mathbf{x}, \mathbf{x} + \mathbf{e}_i - \mathbf{e}_j) = \frac{x_j}{N} \left( \sum_{k=1}^{d} \frac{x_k}{N} m_{ki} \right), \quad 1 \leq i \neq j \leq d;$$

(4.3)
$$K(\mathbf{x}, \mathbf{x}) = 1 - \sum_{i \neq j} K(\mathbf{x}, \mathbf{x} + \mathbf{e}_i - \mathbf{e}_j);$$

$$K(\mathbf{x}, \mathbf{y}) = 0, \quad \text{otherwise.}$$

This model ($d=2$) is due to Moran [39]. Background and references can be found in the text by Ewens [18]. In many applications, the matrix $M = \{m_{ij}\}$ of mutation probabilities takes a special form

$$(4.4) \qquad M = (1-m)I + m \begin{pmatrix} p_1 & \cdots & p_d \\ & \cdots & \\ p_1 & \cdots & p_d \end{pmatrix},$$

where $0 < m < 1$ and $(p_1, \ldots, p_d)$ is a probability vector. In words, the offspring has probability $m$ to mutate. If mutation happens, the offspring will change to species $i$ with probability $p_i$. Note when $m = 1$, the mutation matrix $M$ has identical rows and the process degenerates into a multivariate Ehrenfest chain model, which is a special case in Section 4.3. Karlin and McGregor [36] gave the Karlin–McGregor spectral representation of the transition density for the continuous-time multivariate Moran model. Their version of multivariate Hahn polynomials are defined iteratively and essentially the same as (2.6).

A natural question is how long it takes such a population to be totally mixed. In the univariate ($d=2$) and continuous-time setting, Donnelly and Rodrigues [16] obtain an upper bound of order $(N \log N)/m$ ($m = m_{12} + m_{21}$ in the above notation) in the separation and total variation distances, when the process starts from $X_0 = 0$. As shown below, for the chi-square distance, order $N/m$ (constant being explicit) steps are necessary and sufficient. Convergence rate in the separation distance for the multivariate case is also available [44] but not presented here.

Under the reparametrization $\alpha_i = \frac{Nmp_i}{1-m}$, $1 \leq i \leq d$, it is easy to check that the Moran process defined by (4.3) and (4.4) is exactly the same as the Pólya level model with $s = 1$. This gives an intuitive explanation why the Moran process has Dirichlet-multinomial distribution as its stationary distribution. We remark that all the sequential Pólya urn models admit an interpretation as a population genetics model. For example, a Pólya level model with parameter $s$ means that for a population of size $N$, at each step, $s$ individuals are sequentially selected to reproduce and then mutate, then $s$ individuals are randomly chosen from the old population (size $N$) to die. The original Moran model can be thought of as a multivariate birth–death type process. This generalization allows for more dynamic change of the population at each generation.

The next proposition gives the convergence rate of the Moran process when initially all individuals are of the same species.

PROPOSITION 4.7. *Let $K$ denote the Moran process specified by (4.3) and (4.4). Then $K$ is reversible with respect to the Dirichlet-multinomial distribution $\mathcal{DM}(\cdot|N, \boldsymbol{\alpha})$ (2.4), where $\alpha_i = \frac{Nmp_i}{1-m} > 0$, $1 \leq i \leq d$, and:*

1. *K has eigenvalue*

$$\beta_n = 1 - \frac{n(|\boldsymbol{\alpha}| + n - 1)}{N(N + |\boldsymbol{\alpha}|)}, \qquad 0 \leq n \leq N,$$

*with multiplicity $\binom{d+n-2}{n}$ and corresponding eigenbasis $\{Q_{\mathbf{n}}\}_{|\mathbf{n}|=n}$, where $Q_{\mathbf{n}}$ are the multivariate Hahn polynomials (2.6). Particularly,*

$$\beta_0 = 1 \qquad \text{with multiplicity 1,}$$

$$\beta_1 = 1 - \frac{|\boldsymbol{\alpha}|}{N(N+|\boldsymbol{\alpha}|)} \qquad \text{with multiplicity } d - 1.$$

2. *Suppose that initially all individuals are of species $i$. Then, for any $c > 0$,*

$$\chi^2_{N\mathbf{e}_i}(l) \leq e^{-c}$$

$$\text{for } l \geq \frac{\log[3(2 \vee |\boldsymbol{\alpha}|)N/(N+|\boldsymbol{\alpha}|)((|\boldsymbol{\alpha}| - \alpha_i)/\alpha_i \vee 1)] + c}{-2\log(1 - |\boldsymbol{\alpha}|/(N(N+|\boldsymbol{\alpha}|)))},$$

$$\chi^2_{N\mathbf{e}_i}(l) \geq \frac{1}{6}e^c$$

$$\text{for } l \leq \frac{\log[3(2 \vee |\boldsymbol{\alpha}|)N(|\boldsymbol{\alpha}| - \alpha_i)/((N+|\boldsymbol{\alpha}|)\alpha_i)] - c}{-2\log(1 - |\boldsymbol{\alpha}|/(N(N+|\boldsymbol{\alpha}|)))}.$$

PROOF. The first assertion is a direct corollary to Theorem 4.5 by setting $s = 1$ in (4.1). For the second assertion, by (4.2) and (2.18),

$$\chi^2_{N\mathbf{e}_i}(l) = \sum_{n=1}^N \left(1 - \frac{n(|\boldsymbol{\alpha}| + n - 1)}{N(N + |\boldsymbol{\alpha}|)}\right)^{2l}$$

$$\times \frac{N_{[n]}}{(N+|\boldsymbol{\alpha}|)_{(n)}} \frac{(|\boldsymbol{\alpha}| + 2n - 1)(|\boldsymbol{\alpha}|)_{(n-1)}(|\boldsymbol{\alpha}| - \alpha_i)_{(n)}}{n!(\alpha_i)_{(n)}}.$$

The ratio of the $(n+1)$th summand to the $n$th is

$$\left(1 - \frac{|\boldsymbol{\alpha}| + 2n}{N(N + |\boldsymbol{\alpha}|) - n(|\boldsymbol{\alpha}| + n - 1)}\right)^{2l}$$

$$\times \frac{N-n}{N+|\boldsymbol{\alpha}|+n} \frac{|\boldsymbol{\alpha}|+2n+1}{|\boldsymbol{\alpha}|+2n-1} \frac{|\boldsymbol{\alpha}|+n-1}{n+1} \frac{|\boldsymbol{\alpha}|-\alpha_i+n}{\alpha_i+n}$$

$$\leq (\beta_1)^{2l} \frac{N}{N+|\boldsymbol{\alpha}|} 3\left(1 \vee \frac{|\boldsymbol{\alpha}|}{2}\right) \left(\frac{|\boldsymbol{\alpha}| - \alpha_i}{\alpha_i} \vee 1\right).$$

With $l \geq \frac{\log[3(2\vee|\boldsymbol{\alpha}|)N/(N+|\boldsymbol{\alpha}|)((|\boldsymbol{\alpha}|-\alpha_i)/\alpha_i\vee 1)]+c}{-2\log\beta_1}$, this ratio is bounded by $e^{-c}/2 \leq 1/2$ and hence

$$\chi^2_{N\mathbf{e}_i}(l) \leq \frac{N(|\boldsymbol{\alpha}|+1)(|\boldsymbol{\alpha}|-\alpha_i)}{(N+|\boldsymbol{\alpha}|)\alpha_i}\beta_1^{2l}\left(\sum_{n=0}^\infty \frac{1}{2^n}\right)$$

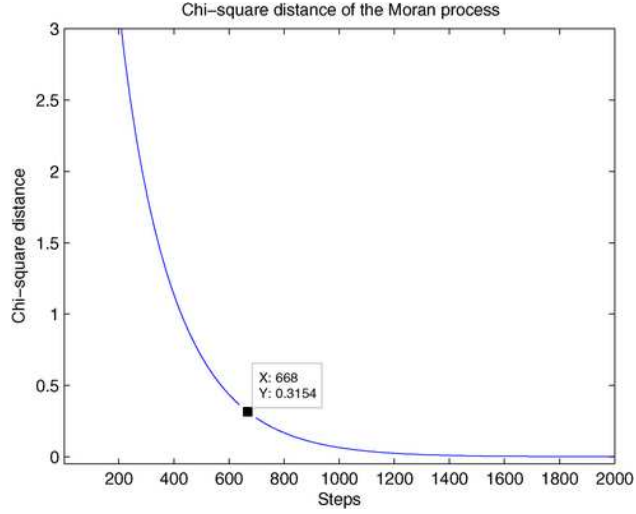

Fig. 3. *Chi-square distance of the Moran model with $N=20$, $d=5$, $\alpha_i = 1/5$, starting from $N\mathbf{e}_i$.*

$$\leq \frac{2(|\boldsymbol{\alpha}|+1)}{3(2\vee|\boldsymbol{\alpha}|)}e^{-c} \leq e^{-c}.$$

With $l \leq \frac{\log[3(2\vee|\boldsymbol{\alpha}|)(|\boldsymbol{\alpha}|-\alpha_i)/\alpha_i N/(N+|\boldsymbol{\alpha}|)]-c}{-2\log\beta_1}$

$$\chi^2_{N\mathbf{e}_i}(l) \geq \frac{N(|\boldsymbol{\alpha}|+1)(|\boldsymbol{\alpha}|-\alpha_i)}{(N+|\boldsymbol{\alpha}|)\alpha_i}\beta_1^{2l} \geq \frac{|\boldsymbol{\alpha}|+1}{3(2\vee|\boldsymbol{\alpha}|)}e^c \geq \frac{1}{6}e^c. \qquad \square$$

REMARK 4.8. In the biological case, $N$ is large and $Nm$ is of constant order. For example, consider the case $m = \frac{1}{N+1}$ and $p_1 = \cdots = p_d = \frac{1}{d}$, or equivalently, $\alpha_1 = \cdots = \alpha_d = \frac{1}{d}$. When initially all individuals are of same type and $N$ is large, $\frac{\log 6(d-1)N/(N+1)}{-2\log(1-1/(N(N+1)))} \approx \frac{\log 6(d-1)}{2}N(N+1)$ steps are necessary and sufficient to drive the chi-square distance low. Figure 3 shows the decrease of the chi-square distance for the Moran process with $N=20$, $d=5$, $\alpha_i = 0.2$. In this case, $\frac{\log 6(d-1)}{2}N(N+1) \approx 668$.

REMARK 4.9. Consider the case $\alpha_1 = \cdots = \alpha_d = 1$. The stationary distribution is uniform on $\mathcal{X}_N^d$. When $N$ is large, $\frac{\log 3d(d-1)N/(N+d)}{-2\log(1-d/(N(N+d)))} \approx \frac{\log 3d(d-1)}{2d} \times N(N+d)$ steps are necessary and sufficient to drive the chi-square distance low. Figure 4 shows the decrease of the chi-square distance for the Moran process with $N=20$, $d=5$, $\alpha_i = 1$. In this case, $\frac{\log 3d(d-1)}{2d}N(N+d) \approx 205$.

4.1.2. *Convergence rate of the local community process by Hubbell.* A description of the local community process by Hubbell is given in the Introduction. We observe that the process is almost the same as the Moran process except that the reproducing individual cannot be the same as the dying one. One-step transition probabilities are

$$K(\mathbf{x}, \mathbf{x} + \mathbf{e}_i - \mathbf{e}_j) = \frac{x_j}{N}\left(\frac{\sum_{k \neq j} x_k m_{ki} + (x_j - 1)m_{ji}}{N-1}\right),$$

(4.5) $$1 \leq i \neq j \leq d;$$

$$K(\mathbf{x}, \mathbf{x}) = 1 - \sum_{i \neq j} K(\mathbf{x}, \mathbf{x} + \mathbf{e}_i - \mathbf{e}_j);$$

$$K(\mathbf{x}, \mathbf{y}) = 0 \quad \text{otherwise},$$

with the matrix of mutation probabilities same as (4.4). This is essentially the process prescribed in Hubbell's book [27], page 86, and simulated by McGill [38]. Again under parametrization $\alpha_i = \frac{(N-1)mp_i}{1-m}$, we find that the local community process by Hubbell is the same as the Pólya down-up model with $s = 1$. The following proposition is similar to Proposition 4.7 and the proof is omitted.

PROPOSITION 4.10. *Let $K$ denote the local community process by Hubbell specified by (4.5) and (4.4). Then $K$ is reversible with respect to the Dirichlet-multinomial distribution $\mathcal{DM}(\cdot|N, \boldsymbol{\alpha})$ (2.4), where $\alpha_i = \frac{(N-1)mp_i}{1-m} > 0$, $1 \leq i \leq d$, and:*

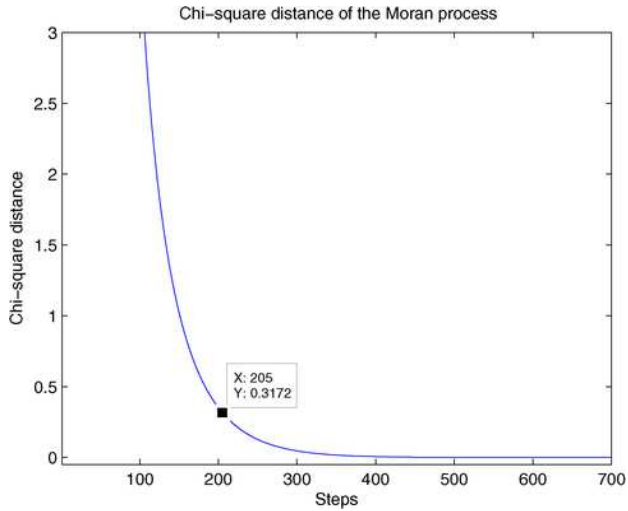

FIG. 4. *Chi-square distance of the Moran model with $N = 20$, $d = 5$, $\alpha_i = 1$, starting from $N\mathbf{e}_i$.*

1. *K has eigenvalue*

$$\beta_n = \frac{(N-n)(N+|\boldsymbol{\alpha}|+n-1)}{N(N+|\boldsymbol{\alpha}|-1)}$$

$$= 1 - \frac{n(n+|\boldsymbol{\alpha}|-1)}{N(N+|\boldsymbol{\alpha}|-1)}, \qquad 0 \leq n \leq N,$$

*with multiplicity $\binom{d+n-2}{n}$ and corresponding eigenbasis $\{Q_\mathbf{n}\}_{|\mathbf{n}|=n}$, where $Q_\mathbf{n}$ are the multivariate Hahn polynomials (2.6). Particularly,*

$$\beta_0 = 1 \qquad \text{with multiplicity } 1,$$

$$\beta_1 = \frac{(N-1)(N+|\boldsymbol{\alpha}|)}{N(N+|\boldsymbol{\alpha}|-1)} = 1 - \frac{|\boldsymbol{\alpha}|}{N(N+|\boldsymbol{\alpha}|-1)} \qquad \text{with multiplicity } d-1.$$

2. *Suppose that initially all individuals are of species $i$. Then for any $c > 0$,*

$$\chi^2_{N\mathbf{e}_i}(l) \leq e^{-c} \qquad \text{for } l \geq \frac{\log[3(2 \vee |\boldsymbol{\alpha}|)N/(N+|\boldsymbol{\alpha}|)((|\boldsymbol{\alpha}|-\alpha_i)/\alpha_i \vee 1)] + c}{-2\log(1-|\boldsymbol{\alpha}|/(N(N+|\boldsymbol{\alpha}|-1)))},$$

$$\chi^2_{N\mathbf{e}_i}(l) \geq \frac{1}{6}e^c \qquad \text{for } l \leq \frac{\log[3(2 \vee |\boldsymbol{\alpha}|)N(|\boldsymbol{\alpha}|-\alpha_i)/((N+|\boldsymbol{\alpha}|)\alpha_i)] - c}{-2\log(1-|\boldsymbol{\alpha}|/(N(N+|\boldsymbol{\alpha}|-1)))}.$$

REMARK 4.11. In the ecological case, $N$ is large and $Nm$ is of constant order. For example, consider the case $m = \frac{1}{N}$ and $p_1 = \cdots = p_d = \frac{1}{d}$, or equivalently, $\alpha_1 = \cdots = \alpha_d = \frac{1}{d}$. When initially all individuals are of the same species and $N$ is large, $\frac{\log((6(d-1)N)/(N+1))}{-2\log(1-1/N^2)} \approx \frac{\log 6(d-1)}{2}N^2$ steps are necessary and sufficient to drive the chi-square distance low.

We have seen a small simulation example in the Introduction. Let us look at another concrete example in [38] where the computation time of the simulation is prohibitive. To sample the stationary distribution of the local community of size $N = 20{,}000$ and migration probability $m = 0.1$ (corresponding to the famous Barro Colorado Island dataset in ecology), McGill first simulates a metacommunity with population size $N_M$ and speciation probability $s$. Given a fixed metacommunity configuration, Proposition 4.10 tells us how many steps need to be run for the local community to reach equilibrium. For example, for a metacommunity configuration with $d = 300$ equally abundant species, 1.44 million steps are necessary and sufficient for the local community to reach equilibrium. This coincides with McGill's empirical findings ([38], Figure 1). If we assume that the tree mortality rate is 1% per year, then this translates into 7200 years for the Barro Colorado Island to reach equilibrium.

4.1.3. *Convergence rate of the Dirichlet-multinomial Gibbs sampler.* This is the multivariate generalization of the canonical beta-binomial Gibbs sampler. The remarkable paper [11] gives explicit diagonalization and sharp convergence rates of Gibbs samplers for six exponential families. One of their motivating examples is how fast the classical beta-binomial Gibbs sampler converges ([11], Proposition 1). We give analogous results for the Dirichlet-multinomial Gibbs sampler here.

Consider the two-component Gibbs sampler with Dirichlet prior $\pi(\mathbf{p}) \sim \mathcal{D}(\cdot|\boldsymbol{\alpha})$ and multinomial likelihood $f(\mathbf{x}|\mathbf{p}) \sim \mathcal{M}(\cdot|N, \mathbf{p})$. The posterior distribution is again Dirichlet, that is, $\pi(\mathbf{p}|\mathbf{x}) \sim \mathcal{D}(\cdot|\mathbf{x} + \boldsymbol{\alpha})$. The Gibbs sampler iterates the following two steps:

- From $\mathbf{x}$, draw $\mathbf{p}$ from $\mathcal{D}(\cdot|\mathbf{x} + \boldsymbol{\alpha})$;
- From $\mathbf{p}$, draw $\mathbf{y}$ from $\mathcal{M}(\cdot|N, \mathbf{p})$.

The marginal $\mathbf{x}$-chain of the Dirichlet-multinomial Gibbs sampler forms a Markov chain with state space $\mathcal{X}_N^d$ and transition probabilities

$$
\begin{aligned}
K(\mathbf{x}, \mathbf{y}) &= \int_{\Delta^d} \mathcal{M}(\mathbf{y}|N, \mathbf{p}) \mathcal{D}(\mathbf{p}|\mathbf{x} + \alpha) \, d\mathbf{p} \\
&= \binom{N}{\mathbf{y}} \frac{\prod_{i=1}^d (x_i + \alpha_i)_{(y_i)}}{(N + |\boldsymbol{\alpha}|)_{(N)}}, \qquad \mathbf{x}, \mathbf{y} \in \mathcal{X}_N^d.
\end{aligned}
$$
(4.6)

We find that the marginal $\mathbf{x}$-chain (4.6) is actually the Pólya level model with $s = N$. Again the proof of the following proposition is analogous to Proposition 4.7 and is omitted.

PROPOSITION 4.12. *Let $K$ denote the marginal $\mathbf{x}$-chain (4.6) of the Dirichlet-multinomial Gibbs sampler. Then $K$ is reversible with respect to the Dirichlet-multinomial distribution $\mathcal{DM}(\cdot|N, \boldsymbol{\alpha})$ (2.4), and:*

1. *$K$ has eigenvalue*

$$\beta_n = \frac{N_{[n]}}{(N + |\boldsymbol{\alpha}|)_{(n)}}, \qquad 0 \leq n \leq N,$$

*with multiplicity $\binom{d+n-2}{n}$ and corresponding eigenbasis $\{Q_{\mathbf{n}}\}_{|\mathbf{n}|=n}$, where $Q_{\mathbf{n}}$ are the multivariate Hahn polynomials (2.6). Particularly,*

$$\beta_0 = 1 \qquad \text{with multiplicity } 1,$$

$$\beta_1 = \frac{N}{N + |\boldsymbol{\alpha}|} \qquad \text{with multiplicity } d - 1.$$

2. *Suppose that the starting state is $N\mathbf{e}_i$. Then for any $c > 0$,*

$$\chi^2_{N\mathbf{e}_i}(l) \leq e^{-c} \qquad for\ l \geq -\frac{1}{2}\left(\frac{\log[3(2 \vee |\boldsymbol{\alpha}|)((|\boldsymbol{\alpha}| - \alpha_i)/\alpha_i \vee 1)] + c}{\log(N/(N + |\boldsymbol{\alpha}|))} + 1\right),$$

$$\chi^2_{N\mathbf{e}_i}(l) \geq \frac{1}{6}e^c \qquad for\ l \leq -\frac{1}{2}\left(\frac{\log[3(2 \vee |\boldsymbol{\alpha}|)(|\boldsymbol{\alpha}| - \alpha_i)/\alpha_i] - c}{\log(N/(N + |\boldsymbol{\alpha}|))} + 1\right).$$

Note that the joint chain $\{(\mathbf{X}_t, \mathbf{P}_t)\}_{t \geq 0}$ is not reversible. Complete analysis of the joint chain relies on a singular value decomposition presented in [11].

REMARK 4.13. Consider the case $\alpha_1 = \cdots = \alpha_d = \frac{1}{d}$. When $N$ is large, $\frac{\log 6(d-1)}{-2\log(N/(N+1))} \approx \frac{\log 6(d-1)}{2}(N+1)$ steps are necessary and sufficient to drive the chi-square distance low. Figure 5 shows the decrease of the chi-square distance for the Dirichlet-multinomial Gibbs sampler with $N = 20$, $d = 5$, $\alpha_i = 0.2$. In this case, $\frac{\log 6(d-1)}{2}(N+1) \approx 33$.

REMARK 4.14. Consider the case $\alpha_1 = \cdots = \alpha_d = 1$ (a uniform sampler on $\mathcal{X}_N^d$). When $N$ is large, $\frac{\log 3d(d-1)}{-2\log(N/(N+d))} \approx \frac{\log 3d(d-1)}{2d}(N+d)$ steps are necessary and sufficient to drive the chi-square distance low. Figure 6 shows the decrease of the chi-square distance for the Dirichlet-multinomial Gibbs sampler with $N = 20$, $d = 5$, $\alpha_i = 1$. In this case, $\frac{\log 3d(d-1)}{2d}(N+d) \approx 10$.

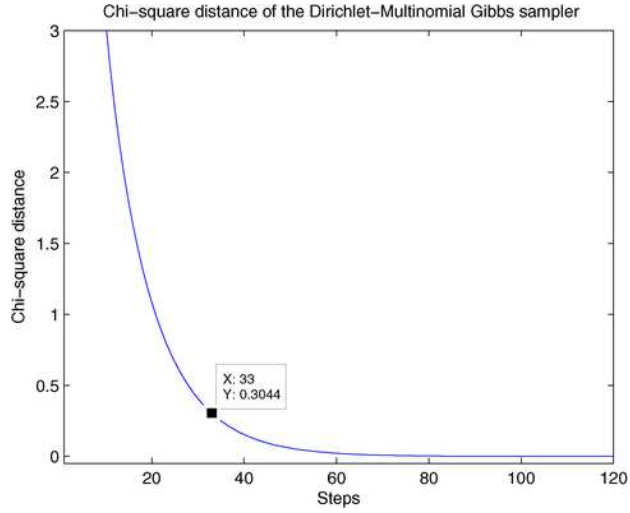

FIG. 5. *Chi-square distance of the Dirichlet-multinomial Gibbs sampler with $N = 20$, $d = 5$, $\alpha_i = 1/5$, starting from $N\mathbf{e}_i$.*

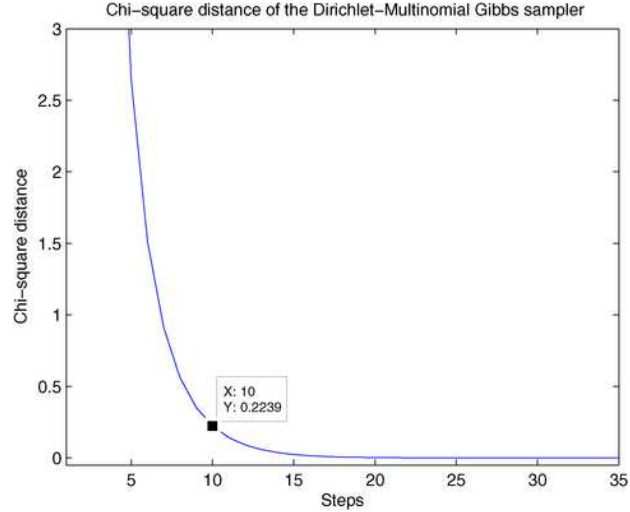

Fig. 6. *Chi-square distance of the Dirichlet-multinomial Gibbs sampler with $N = 20$, $d = 5$, $\alpha_i = 1$, starting from $N\mathbf{e}_i$.*

4.2. *Generalized Bernoulli–Laplace models.* In the classical Bernoulli–Laplace urn model (see, e.g., [31], Section 4.8.1), there are two urns containing $2N$ balls. Initially the left urn contains $N$ red balls; the right urn contains $N$ black balls. At each step, we first randomly choose one ball from the left urn and put it into the right urn. Then we randomly choose one ball from the right urn (it contains $N+1$ balls now) and put it into the left urn. If we track the number of red balls in the left urn, it forms a Markov chain with (univariate) hypergeometric distribution as stationary distribution. In this section, we study generalizations of this model in two directions. First, we allow balls to have more than two colors. Second, we allow more dynamic changes at each step. In analogy to the Pólya urn models, we define three classes of Bernoulli–Laplace models.

There is a batch of balls with composition $\boldsymbol{\ell} = (l_1, \ldots, l_d)$, that is, $l_i$ balls of color $i$ for $1 \leq i \leq d$, distributed in two urns. The left urn contains $N < |\boldsymbol{\ell}|$ balls; the right urn contains $|\boldsymbol{\ell}| - N$ balls. Recall that $|\boldsymbol{\ell}| = \sum_{i=1}^{d} \ell_i$ is the total number of balls.

DEFINITION 4.15. Bernoulli–Laplace level model: $s$ is a parameter satisfying $0 \leq s \leq \min\{N, |\boldsymbol{\ell}| - N\}$. At each step, we randomly choose $s$ balls from each urn and then switch them. Let $\mathbf{X}_t = (X_{t1}, \ldots, X_{td})$ be the composition of the balls in the left urn after $t$ steps. The process $\{\mathbf{X}_t\}_{t \geq 0}$ forms a multivariate Markov chain on $\mathcal{X}^d_{N, \boldsymbol{\ell}}$ and is called a Bernoulli–Laplace level model with parameters $\boldsymbol{\ell}$, $N$ and $s$.

DEFINITION 4.16. Bernoulli–Laplace down-up model: $s$ is a parameter satisfying $0 \leq s \leq N$. At each step, first we randomly choose $s$ balls from the left urn and put them into the right urn. Then we randomly choose $s$ balls from the right urn (it contains $|\boldsymbol{\ell}| - N + s$ balls now) and put them into the left urn. Let $\mathbf{X}_t = (X_{t1}, \ldots, X_{td})$ be the composition of balls in the left urn after $t$ steps. The process $\{\mathbf{X}_t\}_{t \geq 0}$ forms a multivariate Markov chain on $\mathcal{X}^d_{N,\boldsymbol{\ell}}$ and is called a Bernoulli–Laplace down-up model with parameters $\boldsymbol{\ell}$, $N$ and $s$.

DEFINITION 4.17. Bernoulli–Laplace up-down model: $s$ is a parameter satisfying $0 \leq s \leq N$. At each step, first we randomly choose $s$ balls from the right urn and put them into the left urn. Then we randomly choose $s$ balls from the left urn (it contains $N + s$ balls now) and put them into the right urn. Let $\mathbf{X}_t = (X_{t1}, \ldots, X_{td})$ be the composition of balls in the left urn after $t$ steps. The process $\{\mathbf{X}_t\}_{t \geq 0}$ forms a multivariate Markov chain on $\mathcal{X}^d_{N,\boldsymbol{\ell}}$ and is called a Bernoulli–Laplace up-down model with parameters $\boldsymbol{\ell}$, $N$ and $s$.

The special case $d = 2$, $s = 1$ of the down-up model dates back to Bernoulli and Laplace who introduced this model to study diffusion of particles between two containers. More details and historical background can be found in [20, 31]. Both [6] and [15] study convergence rates of the Bernoulli–Laplace level model with $d = 2$, $s = 1$ and contain interesting connections to real world problems.

This process can be lifted to a random walk on the space of all $N$-subsets of $|\boldsymbol{\ell}|$ objects. For example, in the $s = 1$ case of the level model, at each step, randomly pick one element from the current set, one from the complement set and switch them. This is a nearest-neighbor random walk under the metric $d(x, y) = N - |x \cap y|$. The stationary distribution is uniform over all $N$-subsets. Therefore the Bernoulli–Laplace process has hypergeometric distribution as stationary distribution. See [6] for detailed analysis of the lifted chain. This point of view gives the following result.

LEMMA 4.18. *All three classes of Bernoulli–Laplace models are reversible with respect to the multivariate hypergeometric distribution $\mathcal{H}(\cdot|N, \boldsymbol{\ell})$ (2.5).*

The Bernoulli–Laplace up-down and down-up models have alternative interpretations as the marginal chains of a multivariate hypergeometric walk. The univariate hypergeometric walk was originally studied in [14]. Consider the following Gibbs sampler. Let $\boldsymbol{\ell} = (l_1, \ldots, l_d)$ be a vector of counts and $0 < N \leq N + s \leq |\boldsymbol{\ell}|$. The prior distribution is hypergeometric:

$$\pi(\boldsymbol{\theta}) \sim \mathcal{H}(\cdot|N + s, \boldsymbol{\ell}),$$

that is, $\boldsymbol{\theta} = (\theta_1, \ldots, \theta_d)$ is the composition of a random sample of size $N + s$ from a pool of $|\boldsymbol{\ell}|$ balls with composition $\boldsymbol{\ell}$. The likelihood given parameter $\boldsymbol{\theta}$ is

$$f(\mathbf{x}|\boldsymbol{\theta}) \sim \mathcal{H}(\cdot|N, \boldsymbol{\theta}),$$

that is, $\mathbf{X} = (X_1, \ldots, X_d)$ is the composition of a random sample of size $N$ from a pool of $|\boldsymbol{\theta}| = N + s$ balls with composition $\boldsymbol{\theta}$. $\pi(\boldsymbol{\theta})$ is the conjugate prior for the hypergeometric distribution and the posterior distribution is still hypergeometric:

$$\pi(\boldsymbol{\theta}|\mathbf{x}) \sim \mathbf{x} + \mathcal{H}(\cdot|s, \boldsymbol{\ell} - \mathbf{x}),$$

that is, take a random sample of size $s$ from a pool of $|\boldsymbol{\ell}| - N$ balls with composition $\boldsymbol{\ell} - \mathbf{x}$ and set $\theta_i$ to be the count of balls of color $i$ plus $x_i$. The marginal $\mathbf{x}$-chain has transition kernel

$$K(\mathbf{x}, \mathbf{y}) = \sum_{\boldsymbol{\theta}} \pi(\boldsymbol{\theta}|\mathbf{x}) f(\mathbf{y}|\boldsymbol{\theta})$$

and is the same as a Bernoulli–Laplace up-down model with parameters $\boldsymbol{\ell}, N, s$. The marginal $\boldsymbol{\theta}$-chain has transition kernel

$$K(\boldsymbol{\theta}, \boldsymbol{\theta}') = \sum_{\mathbf{x} \in \mathcal{X}_{N,\boldsymbol{\theta}}^d} f(\mathbf{x}|\boldsymbol{\theta}) \pi(\boldsymbol{\theta}'|\mathbf{x})$$

and is the same as a Bernoulli–Laplace down-up model with parameters $\boldsymbol{\ell}, N + s, N$.

In analogy to the Pólya urn models, these models share the same polynomial eigenfunctions which are the multivariate Hahn polynomials for the hypergeometric distribution. The following theorem is analogous to Theorem 4.5 and the proof is omitted.

THEOREM 4.19. 1. *The Bernoulli–Laplace level model has eigenvalue*

$$\beta_n = \sum_{k=0}^{n} \binom{n}{k} \frac{(N-s)_{[k]} s_{[n-k]}}{N_{[k]}(|\boldsymbol{\ell}| - N)_{[n-k]}}$$
$$= \sum_{k=0}^{n} \binom{n}{k} \frac{(N-s)_{[n-k]} s_{[k]}}{N_{[n-k]}(|\boldsymbol{\ell}| - N)_{[k]}}, \qquad 0 \le n \le N,$$

*with multiplicity $|\mathcal{X}_{n,\boldsymbol{\ell}}^d|$ and corresponding eigenbasis $\{Q_{\mathbf{n}}(\mathbf{x}; N, -\boldsymbol{\ell})\}_{|\mathbf{n}|=n}$, where $Q_{\mathbf{n}}$ are the multivariate Hahn polynomials for the hypergeometric distribution as defined in Proposition 2.3.*

2. *The Bernoulli–Laplace down-up model has eigenvalue*

$$\beta_n = \frac{(N-s)_{[n]}(|\boldsymbol{\ell}| - N)_{[n]}}{N_{[n]}(|\boldsymbol{\ell}| - N + s)_{[n]}}, \qquad 0 \le n \le N,$$

with multiplicity $|\mathcal{X}_{n,\boldsymbol{\ell}}^d|$ and corresponding eigenbasis $\{Q_{\mathbf{n}}(\mathbf{x}; N, -\boldsymbol{\ell})\}_{|\mathbf{n}|=n}$, where $Q_{\mathbf{n}}$ are the multivariate Hahn polynomials for the hypergeometric distribution as defined in Proposition 2.3.

3. *The Bernoulli–Laplace up-down model has eigenvalue*

$$\beta_n = \frac{N_{[n]}(|\boldsymbol{\ell}| - N - s)_{[n]}}{(N+s)_{[n]}(|\boldsymbol{\ell}| - N)_{[n]}}, \qquad 0 \leq n \leq N,$$

with multiplicity $|\mathcal{X}_{n,\boldsymbol{\ell}}^d|$ and corresponding eigenbasis $\{Q_{\mathbf{n}}(\mathbf{x}; N, -\boldsymbol{\ell})\}_{|\mathbf{n}|=n}$, where $Q_{\mathbf{n}}$ are the multivariate Hahn polynomials for the hypergeometric distribution as defined in Proposition 2.3.

4. *For all three classes of Bernoulli–Laplace models, the chi-square distance between $K^l(\mathbf{x}, \cdot)$ and stationary distribution is*

$$\chi_{\mathbf{x}}^2(l) = \sum_{n=1}^{N} \beta_n^{2l} h_n(\mathbf{x}, \mathbf{x}),$$

where $\beta_n$ are the eigenvalues for the corresponding process and $h_n$ are the kernel polynomials for the hypergeometric distribution, that is, (2.16) with $\alpha_i$ replaced by $-l_i$.

REMARK 4.20. For all three classes of the Bernoulli–Laplace urn models, the eigenvalues depend on parameters $\boldsymbol{\ell} = (l_1, \ldots, l_d)$ only through $|\boldsymbol{\ell}|$. Formulas for the eigenfunctions require complete knowledge of $\boldsymbol{\ell}$.

4.2.1. *Convergence rate of the Bernoulli–Laplace down-up models.* We specialize to the case in which at beginning the left urn contains balls of the same color.

PROPOSITION 4.21. *For $N \leq l_i$ and any $c > 0$, the chi-square distance between the Bernoulli–Laplace down-up model and the stationary distribution, starting from $N\mathbf{e}_i$, satisfies*

$$\chi_{N\mathbf{e}_i}^2(l) \leq 2e^{-c} \qquad \text{for } l \geq \frac{\log[|\boldsymbol{\ell}|N(|\boldsymbol{\ell}| - l_i)/((|\boldsymbol{\ell}| - N)l_i)] + c}{-2\log[(N-s)(|\boldsymbol{\ell}| - N)/(N(|\boldsymbol{\ell}| - N + 1))]},$$

$$\chi_{N\mathbf{e}_i}^2(l) \geq \frac{1}{2}e^c \qquad \text{for } l \leq \frac{\log[|\boldsymbol{\ell}|N(|\boldsymbol{\ell}| - l_i)/((|\boldsymbol{\ell}| - N)l_i)] - c}{-2\log[(N-s)(|\boldsymbol{\ell}| - N)/(N(|\boldsymbol{\ell}| - N + 1))]}.$$

The proof is similar to that for the sequential Pólya urn models and is omitted. This result is useful when $s \ll N$. In the extreme case $s = N$, the chain achieves stationarity after one-step.

REMARK 4.22. Consider the multivariate version of the classical Bernoulli–Laplace model ([31], Section 4.8.1), where $s = 1$. Suppose the two urns have equal sizes $N = |\ell|/2$ and initially the left urn contains $N$ balls of the same color and the right urn contains $N$ balls of colors different from those in the left urn ($\mathbf{X}_0 = N\mathbf{e}_i$, $l_i = N$). Then, for $N$ large, $\frac{\log 2N}{-2\log(1-2/(N+1))} \approx \frac{(N+1)\log 2N}{4}$ steps are necessary and sufficient to drive the chi-square distance low. Figure 7 shows the decrease of the chi-square distance for the Bernoulli–Laplace process with $N = 20$, $d = 5$, $|\ell| = 2N$, $l_i = N$. In this case, $\frac{(N+1)\log 2N}{4} \approx 19$.

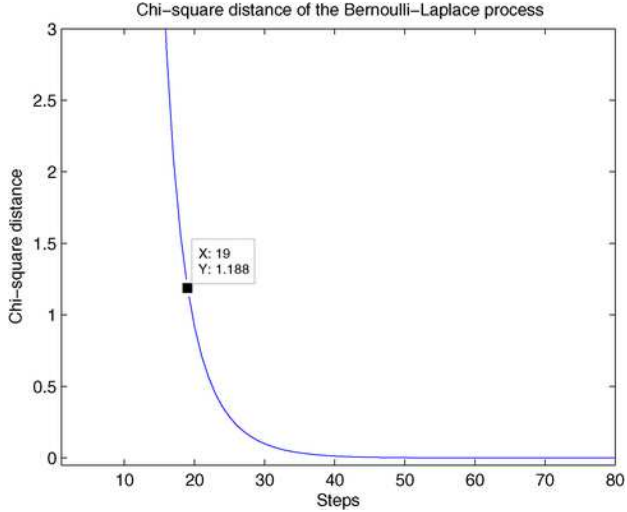

Fig. 7. *Chi-square distance of the Bernoulli–Laplace process with $N = 20$, $d = 5$, $|\ell| = 2N, l_i = N$, starting from $N\mathbf{e}_i$.*

4.3. *A generalized Ehrenfest urn model.* In the classical Ehrenfest model (see, e.g., [20, 34]), a certain number of balls are shuttled between two urns. At each step one ball is randomly chosen and shifted to the other urn. The Ehrenfest chain tracks the number of balls in one of the urns. Reference [35] contains a discussion of the Ehrenfest urn model with $d > 2$ urns. The discrete-time analog of their continuous-time Markov chain randomly chooses a single ball at each step and redistributes it to the $d$ urns according to a probability vector $\mathbf{p} = (p_1, \ldots, p_d)$. The multivariate Ehrenfest chain tracks the counts in each urn.

In this section, we generalize even further by selecting more balls at each step. Specifically, there are $N$ indistinguishable balls distributed in $d$ urns. $s \in \{0, 1, \ldots, N\}$ is a parameter. At each step, we randomly choose $s$ balls from the total of $N$ balls and redistribute each of them independently according to the same probability vector $\mathbf{p}$. Let $X_{ti}$ be the number of balls in the $i$th urn at time $t$. Then $\{\mathbf{X}_t = (X_{t1}, \ldots, X_{td})\}_{t \geq 0}$ forms a multivariate Markov chain on $\mathcal{X}_N^d$.

It is interesting to observe that this generalized Ehrenfest model has an alternative interpretation as the marginal chain of a Gibbs sampler. Consider a multinomial sampling model (sample size $s$) with a preexisting sample of size $(N - s)$. The likelihood of the data given parameter $\boldsymbol{\theta} \in \mathcal{X}_{N-s}^d$ is

$$f(\mathbf{x}|\boldsymbol{\theta}) \sim \boldsymbol{\theta} + \mathcal{M}(\cdot|s, \mathbf{p}).$$

If we specify the prior distribution for the parameter $\boldsymbol{\theta}$ as

$$\pi(\boldsymbol{\theta}) \sim \mathcal{M}(\cdot|N - s, \mathbf{p}),$$

then by the Bayes formula, the posterior distribution of $\theta$ given the data has probability mass function

$$\pi(\boldsymbol{\theta}|\mathbf{x}) = \frac{\binom{N-s}{\boldsymbol{\theta}}\binom{s}{\mathbf{x}-\boldsymbol{\theta}}}{\binom{N}{\boldsymbol{\theta}}}.$$

Consider the Gibbs sampling procedure to sample from the joint distribution of $(\mathbf{X}, \boldsymbol{\theta})$. The Gibbs sampler iterates the following two steps:

- From $\mathbf{x}$, draw $\boldsymbol{\theta}$ from the posterior $\pi(\cdot|\mathbf{x})$.
- From $\boldsymbol{\theta}$, draw $\mathbf{y}$ from the likelihood function $f(\cdot|\boldsymbol{\theta})$.

The marginal $\mathbf{x}$-chain has transition probabilities

$$K(\mathbf{x}, \mathbf{y}) = \sum_{\theta \in \mathcal{X}_{N-s}^d} \pi(\boldsymbol{\theta}|\mathbf{x}) f(\mathbf{y}|\boldsymbol{\theta}), \qquad \mathbf{x}, \mathbf{y} \in \mathcal{X}_N^d,$$

and is the same as the generalized Ehrenfest urn model described previously. Choosing $\boldsymbol{\theta}$ from $\pi(\cdot|\mathbf{x})$ corresponds to choosing $N - s$ balls which will not be redistributed (and hence $s$ balls which will be redistributed). Choosing $\mathbf{y}$ from $f(\cdot|\boldsymbol{\theta})$ corresponds to redistributing $s$ balls independently according to the same probability vector $\mathbf{p}$.

It is well known that the marginal $\mathbf{x}$-chain of a Gibbs sampling Markov chain is reversible, with the marginal of the joint distribution of $(\mathbf{X}, \boldsymbol{\theta})$ as its stationary distribution.

LEMMA 4.23. *The generalized Ehrenfest urn model is reversible with respect to the multinomial distribution $\mathcal{M}(\cdot|N, \mathbf{p})$.*

The following result is again an application of Lemma 3.2.

THEOREM 4.24. 1. *The generalized Ehrenfest urn model has eigenvalue*

$$\beta_n = \frac{(N-s)_{[n]}}{N_{[n]}}, \qquad 0 \leq n \leq N,$$

with multiplicity $\binom{d+n-2}{n}$ and corresponding eigenbasis $\{K_{\mathbf{n}}\}_{|\mathbf{n}|=n}$, where $K_{\mathbf{n}}$ are the multivariate Krawtchouk polynomials (2.10) and $|\mathbf{n}| = \sum_{i=1}^{d} n_i$. In particular,

$$\beta_0 = 1 \qquad \text{with multiplicity } 1,$$
$$\beta_1 = 1 - \frac{s}{N} \qquad \text{with multiplicity } d-1.$$

2. If the process starts from $\mathbf{x}$, then the chi-square distance after $l$ steps is

(4.7) $$\chi^2_{\mathbf{x}}(l) = \sum_{n=1}^{N} \beta_n^{2l} h_n(\mathbf{x}, \mathbf{x}),$$

where $h_n$ are the kernel polynomials (2.22) for the multinomial distribution $\mathcal{M}(\cdot | N, \mathbf{p})$.

We observe that, in terms of convergence to stationarity, the worst-case initial configuration is one of the $d$ configurations where all the balls are in a single urn.

COROLLARY 4.25. *Suppose that initially all balls are in the $i$th urn. Then, for any $c > 0$, the chi-square distance of the generalized Ehrenfest chain satisfies*

$$\chi^2_{N\mathbf{e}_i}(l) \leq e^{e^{-c}} - 1 \qquad \text{for } l \geq \frac{\log(N(1-p_i)/p_i) + c}{-2\log(1-s/N)},$$

$$\chi^2_{N\mathbf{e}_i}(l) \geq e^c \qquad \text{for } l \leq \frac{\log(N(1-p_i)/p_i) - c}{-2\log(1-s/N)}.$$

PROOF. Note $\beta_n = 0$ for $n > N - s$. By (4.7) and (2.23),

$$\chi^2_{N\mathbf{e}_i}(l) = \sum_{n=1}^{N-s} \left(\frac{(N-s)_{[n]}}{N_{[n]}}\right)^{2l} \binom{N}{n} \left(\frac{1-p_i}{p_i}\right)^n.$$

The inequality $\frac{(N-s)_{[n]}}{N_{[n]}} \leq (1 - \frac{s}{N})^n$ implies that

$$N\left(\frac{1-p_i}{p_i}\right)\left(1 - \frac{s}{N}\right)^{2l} \leq \chi^2_{N\mathbf{e}_i}(l) \leq \left(1 + \frac{1-p_i}{p_i}\left(1 - \frac{s}{N}\right)^{2l}\right)^N - 1.$$

Substituting $l = \frac{\log(N(1-p_i)/p_i) + c}{-2\log(1-s/N)}$, we get

$$e^{-c} \leq \chi^2_{N\mathbf{e}_i}(l) \leq \left(1 + \frac{e^{-c}}{N}\right)^N - 1 \leq e^{e^{-c}} - 1.$$

And substituting $l = \frac{\log(N(1-p_i)/p_i) - c}{-2\log(1-s/N)}$, we get

$$\chi^2_{N\mathbf{e}_i}(l) \geq e^c.$$

□

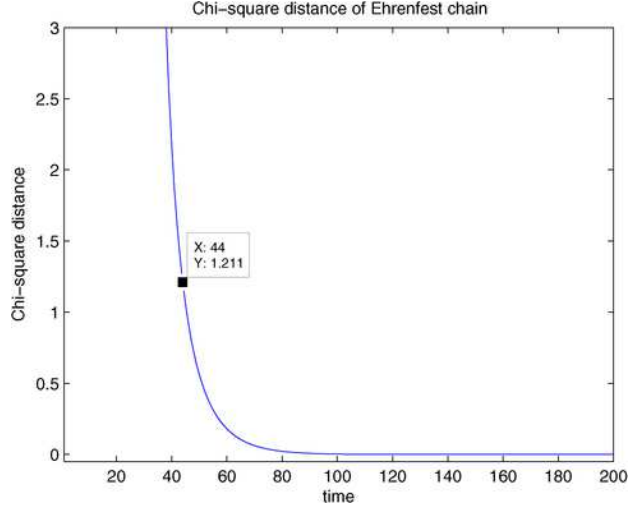

Fig. 8. *Chi-square distance of the Ehrenfest urn model with $N = 20$, $s = 1$, $d = 5$, $p_i = 1/5$, starting from $N\mathbf{e}_i$.*

This analysis is particularly useful when $s$ is small compared to $N$. In the extreme case $s = N$, the chain achieves stationarity in just one-step.

REMARK 4.26. Consider the discrete-time version of the multivariate Ehrenfest urn model in [35] where $s = 1$. For $N$ large, $\frac{\log(N(1-p_i)/p_i)}{-2\log(1-1/N)} \approx \frac{N \log(N(1-p_i)/p_i)}{2}$ steps are necessary and sufficient to drive the chi-square distance low. A probabilistic analysis of this chain is given in [10]. Figure 8 shows the decrease of the chi-square distance for the Ehrenfest urn model with $N = 20$, $s = 1$, $d = 5$, $p_i = 1/5$. In this case, $\frac{N \log(N(1-p_i)/p_i)}{2} \approx 44$.

4.4. *Multivariate normal autoregressive process.* Consider a multivariate normal autoregressive process on $\mathbb{R}^d$ defined by

$$\mathbf{X}_t = A\mathbf{X}_{t-1} + \boldsymbol{\xi}_t, \qquad t \geq 1, \tag{4.8}$$

where $\{\boldsymbol{\xi}_t\}_{t \geq 1}$ are independent and identically distributed $\mathcal{N}(\mathbf{0}, V)$. This process arises in the multigrid Monte Carlo method by Goodman and Sokal [22] and general overrelaxation MCMC algorithms [3, 4, 40]. First we check the stationary distribution of the process and the reversibility criterion.

PROPOSITION 4.27. *The Markov chain (4.8) has a unique stationary distribution $\mathcal{N}(\mathbf{0}, \Sigma)$ if and only if $V = \Sigma - A\Sigma A^T$. Moreover, when $V = \Sigma - A\Sigma A^T$, the Markov chain is reversible if and only if $A\Sigma = \Sigma A^T$.*

PROOF. If $\mathbf{X}_t \sim \mathcal{N}(\mathbf{0}, \Sigma)$, then $\mathbf{X}_{t+1} = A\mathbf{X}_t + \boldsymbol{\xi}_{t+1} \sim \mathcal{N}(\mathbf{0}, V + A\Sigma A^T)$. Therefore $\mathcal{N}(\mathbf{0}, \Sigma)$ is a stationary distribution if and only if $\Sigma = V + A\Sigma A^T$. For the uniqueness, since $V = \Sigma - A\Sigma A^T$ is a positive definite matrix, the spectral radius of $A$ is strictly less than 1. By iterating, $\mathbf{X}_t|\mathbf{X}_0 = \mathbf{x} \sim \mathcal{N}(A^t \mathbf{x}, \sum_{j=0}^{t-1} A^j V (A^T)^j)$. But $A^t \mathbf{x} \to 0$ and

$$\sum_{j=0}^{t-1} A^j V (A^T)^j = \sum_{j=0}^{t-1} [A^j \Sigma (A^T)^j - A^{j+1} \Sigma (A^T)^{j+1}]$$
$$= \Sigma - A^t \Sigma (A^T)^t$$
$$\to \Sigma$$

as $t \to \infty$. Therefore $\mathcal{N}(\mathbf{0}, \Sigma)$ is the unique stationary distribution.

To check the reversibility criterion, note the transition density and stationary density are

$$K(\mathbf{x}, \mathbf{y}) = \frac{1}{\sqrt{(2\pi)^d |V|}} e^{(\mathbf{y} - A\mathbf{x})^T V^{-1}(\mathbf{y} - A\mathbf{x})/2},$$

$$\pi(\mathbf{y}) = \frac{1}{\sqrt{(2\pi)^d |\Sigma|}} e^{-\mathbf{y}^T \Sigma^{-1} \mathbf{y}/2}, \qquad \mathbf{x}, \mathbf{y} \in \mathbb{R}^d.$$

Hence, the Markov chain is reversible if and only if

$$\pi(\mathbf{x}) K(\mathbf{x}, \mathbf{y}) = \pi(\mathbf{y}) K(\mathbf{y}, \mathbf{x}) \quad \text{for all } \mathbf{x}, \mathbf{y}$$
$$\Leftrightarrow \quad (\mathbf{y} - A\mathbf{x})^T V^{-1} (\mathbf{y} - A\mathbf{x}) + \mathbf{x}^T \Sigma^{-1} \mathbf{x}$$
$$\qquad = (\mathbf{x} - A\mathbf{y})^T V^{-1} (\mathbf{x} - A\mathbf{y}) + \mathbf{y}^T \Sigma^{-1} \mathbf{y} \quad \text{for all } \mathbf{x}, \mathbf{y}$$
$$\Leftrightarrow \quad \Sigma^{-1} + A^T V^{-1} A = V^{-1}, \qquad A^T V^{-1} = V^{-1} A$$
$$\Leftrightarrow \quad \Sigma^{-1} + V^{-1} A^2 = V^{-1}, \qquad A^T V^{-1} = V^{-1} A$$
$$\Leftrightarrow \quad V + A^2 \Sigma = \Sigma, \qquad A^T V^{-1} = V^{-1} A$$
$$\Leftrightarrow \quad \Sigma - A\Sigma A^T + A^2 \Sigma = \Sigma, \qquad V A^T = AV$$
$$\Leftrightarrow \quad A\Sigma A^T = A^2 \Sigma, \qquad \Sigma A^T - A\Sigma (A^T)^2 = A\Sigma - A^2 \Sigma A^T$$
$$\Leftrightarrow \quad A\Sigma A^T = A^2 \Sigma, \qquad \Sigma A^T = A\Sigma$$
$$\Leftrightarrow \quad A\Sigma = \Sigma A^T. \qquad \square$$

In the following, we assume that $V = \Sigma - A\Sigma A^T$ and $A\Sigma = \Sigma A^T$. Hence, the Markov chain (4.8) has a unique stationary distribution $\pi \sim \mathcal{N}(\mathbf{0}, \Sigma)$ and is reversible. Let $\lambda_1$ be the largest eigenvalue of $A$. In [22], Goodman and Sokal identify $\lambda_1$ as the rate of decay of the autocorrelation functions of $\mathbf{X}_t$. In [4], Barone and Frigessi identify $\lambda_1$ as the rate of variational norm convergence. Roberts and Sahu [40] show that for all $f \in l^2(\pi)$ and $r > \lambda_1$,

$$\lim_{t \to \infty} \frac{\mathbf{E}_\pi[(f(\mathbf{X}_t) - \int f d\pi)^2]}{r^t} = 0.$$

Our contribution is to use the eigenfunctions of this Markov chain to obtain an exact expression for the chi-square distance from stationarity after $\ell$ steps, for any $\ell \geq 0$ (Proposition 4.30). This leads to nonasymptotic bounds for convergence to stationarity and strengthens earlier results.

The idea is that the original autoregressive process can be transformed into another that can be easily analyzed. The condition $A\Sigma = \Sigma A^T$ implies that $\Sigma^{-1/2} A \Sigma^{1/2}$ is symmetric and thus orthogonally diagonalizable. Let $\Sigma^{-1/2} A \Sigma^{1/2} = PDP^T$ be its eigendecomposition, where $P^T P = I$ and $D$ is the diagonal matrix containing the eigenvalues $1 > |\lambda_1| \geq \cdots \geq |\lambda_d|$ of $A$. We study the transformed Markov chain $\{\mathbf{Z}_t\}$, where

$$\mathbf{Z}_t = P^T \Sigma^{-1/2} \mathbf{X}_t, \qquad t \geq 0.$$

From (4.8),

$$\mathbf{Z}_t = (P^T \Sigma^{-1/2} A \Sigma^{1/2} P) \mathbf{Z}_{t-1} + (P^T \Sigma^{-1/2}) \boldsymbol{\xi}_t, \qquad t \geq 0.$$

Note that

$$P^T \Sigma^{-1/2} A \Sigma^{1/2} P = D,$$
$$\mathrm{Var}(P^T \Sigma^{-1/2} \boldsymbol{\xi}_n) = P^T \Sigma^{-1/2} C \Sigma^{-1/2} P$$
$$= P^T \Sigma^{-1/2} (\Sigma - A \Sigma A^T) \Sigma^{-1/2} P$$
$$= I - D^2.$$

It follows that

(4.9) $$\mathbf{Z}_t = D\mathbf{Z}_{t-1} + \boldsymbol{\xi}'_t, \qquad t \geq 0,$$

where $\boldsymbol{\xi}'_t$ are independent and identically distributed $\mathcal{N}(\mathbf{0}, I - D^2)$. Since $D$ is a diagonal matrix with entries $\lambda_1, \lambda_2, \ldots, \lambda_d$ and components of $\xi'_t$ are independent, *all components of the $\mathbf{Z}_t$ chain proceed as independent univariate normal autoregressive processes*, that is, for $1 \leq i \leq d$,

(4.10) $$Z_{t,i} = \lambda_i Z_{t-1,i} + \xi'_{t,i},$$

where $\xi'_{t,i}$, $t \geq 1$, are independent and identically distributed $\mathcal{N}(0, 1 - \lambda_i^2)$.

Univariate normal autoregressive processes are well studied. The $i$th component process (4.10) is reversible with respect to the standard normal distribution $\mathcal{N}(0,1)$ and has eigenvalues $\lambda_i^n, n \geq 0$, with the Hermit polynomials $\{H_n\}_{n \geq 0}$ (2.14) as the corresponding eigenfunctions. This spectral information is easily transferred to that of the product chain due to independence.

LEMMA 4.28. *Let $(K', \pi')$ be the Markov operator and stationary distribution corresponding to the transformed Markov chain $\{\mathbf{Z}_t\}_{t \geq 0}$ in (4.9). Then $K'$ is reversible with respect to $\pi' \sim \mathcal{N}(\mathbf{0}, I)$. Moreover,*

1. *$K'$ has eigenvalues*

$$\beta_{\mathbf{n}} = \prod_{i=1}^{d} \lambda_i^{n_i}, \qquad \mathbf{n} = (n_1, \ldots, n_d) \in \mathbb{N}_0^d,$$

*with corresponding eigenfunctions $H_{\mathbf{n}}(\frac{\mathbf{z}}{\sqrt{2}}) = \prod_{i=1}^{d} H_{n_i}(\frac{z_i}{\sqrt{2}})$, which satisfy the orthogonality relation*

$$\int_{\mathbb{R}^d} H_{\mathbf{n}}\left(\frac{\mathbf{z}}{\sqrt{2}}\right) H_{\mathbf{m}}\left(\frac{\mathbf{z}}{\sqrt{2}}\right) \frac{e^{-\mathbf{z}^T \mathbf{z}/2}}{(\sqrt{2\pi})^d} d\mathbf{z} = 2^{|\mathbf{n}|} \prod_{i=1}^{d} n_i! \delta_{\mathbf{nm}}.$$

2. *The chi-square distance of $K'$, starting from state $\mathbf{z} \in \mathbb{R}^d$, after $l$ steps is*

$$\chi_{\mathbf{z}}^2(l) = \frac{e^{\sum_{i=1}^{d} z_i^2 \lambda_i^{2l}/(1+\lambda_i^{2l})}}{\sqrt{\prod_{i=1}^{d}(1-\lambda_i^{4l})}} - 1.$$

PROOF. Let $K_i$ denote the Markov operator of the $i$th component process. By independence between the component processes

$$KH_{\mathbf{n}}\left(\frac{\mathbf{z}}{\sqrt{2}}\right) = \mathbf{E}_{K(\mathbf{z},\cdot)}\left[H_{\mathbf{n}}\left(\frac{\mathbf{Z}}{\sqrt{2}}\right)\right] = \prod_{i=1}^{d} \mathbf{E}_{K_i(z_i,\cdot)}\left[H_{n_i}\left(\frac{Z_i}{\sqrt{2}}\right)\right]$$

$$= \prod_{i=1}^{d} K_i H_{n_i}\left(\frac{z_i}{\sqrt{2}}\right) = \left(\prod_{i=1}^{d} \gamma_i^{n_i}\right) H_{\mathbf{n}}\left(\frac{\mathbf{z}}{\sqrt{2}}\right).$$

Then the first assertion follows. For calculation of the chi-square distance,

$$\chi_{\mathbf{z}}^2(l) = \sum_{\mathbf{n} \neq \mathbf{0}} \beta_{\mathbf{n}}^{2l} H_{\mathbf{n}}^2\left(\frac{\mathbf{z}}{\sqrt{2}}\right) d_{\mathbf{n}}^{-2} = \sum_{\mathbf{n} \in \mathbb{N}_0^d} \prod_{i=1}^{d} \frac{\lambda_i^{2n_i l} H_{n_i}^2(z_i/\sqrt{2})}{2^{n_i} n_i!} - 1$$

$$= \prod_{i=1}^{d} \sum_{n_i=0}^{\infty} \frac{\lambda_i^{2n_i l} H_{n_i}^2(z_i/\sqrt{2})}{2^{n_i} n_i!} - 1 = \frac{e^{\sum_{i=1}^{d} z_i^2 \lambda_i^{2l}/(1+\lambda_i^{2l})}}{\sqrt{\prod_{i=1}^{d}(1-\lambda_i^{4l})}} - 1.$$

The last equality follows from the multilinear generating function (2.15). □

It is trivial to check that the chi-square distance of the transformed chain is equal to that of the original chain $K$. This implies the following result.

PROPOSITION 4.29. *For the Markov chain (4.8) satisfying $C = \Sigma - A\Sigma A^T$ and $A\Sigma = \Sigma A^T$, the chi-square distance after $l$ steps, starting from the state $\mathbf{x} \in \mathbb{R}^d$, is*

$$\chi^2_{\mathbf{x}}(l) = \frac{e^{\sum_{i=1}^d z_i^2 \lambda_i^{2l}/(1+\lambda_i^{2l})}}{\sqrt{\prod_{i=1}^d (1-\lambda_i^{4l})}} - 1,$$

*where $\mathbf{z} = P^T \Sigma^{-1/2} \mathbf{x}$.*

Clearly, all the eigenvalues $(\lambda_1, \ldots, \lambda_d)$ of $A$ play a role in determining the speed of convergence. However, if one is willing to compromise a little on sharpness of the bounds, we can get a result only involving the largest eigenvalue $\lambda_1$ of $A$.

PROPOSITION 4.30. *For the Markov chain (4.8), when starting from the state $\mathbf{0}$,*

$$\chi^2_{\mathbf{0}}(l) \leq 10 e^{-c} \quad \text{for } l \geq \frac{\log 2 + c}{-4 \log |\lambda_1|}, \ c \geq \log\left(\frac{d}{2}\right);$$

$$\chi^2_{\mathbf{0}}(l) \geq \frac{1}{4} e^c \quad \text{for } l \leq \frac{\log 2 - c}{-4 \log |\lambda_1|}, \ c > 0.$$

PROOF. Note that $\chi^2_{\mathbf{0}}(l) \leq (1 - \lambda_1^{4l})^{-d/2} - 1$. If $l \geq \frac{\log 2 + c}{-4 \log |\lambda_1|}$, then $\lambda_1^{4l} \leq e^{-c}/2 \leq 1/2$. Hence,

$$\chi^2_{\mathbf{0}}(l) \leq (1 - \lambda_1^{4l})^{-d/2} - 1 \leq (1 + e^{-c})^d - 1$$
$$\leq e^{de^{-c}} e^{-c} \leq 10 e^{-c},$$

when $c \geq \log(d/2)$. For the lower bound,

$$\chi^2_{\mathbf{0}}(l) \geq (1 - \lambda_1^{4l})^{-1/2} - 1 \geq \frac{\lambda_1^{4l}}{2}.$$

Hence, if $l \leq \frac{\log 2 - c}{-4 \log |\lambda_1|}$,

$$\chi^2_{\mathbf{0}}(l) \geq \frac{e^c}{4}. \qquad \square$$

4.4.1. *An example from image analysis.* We now borrow an example from Bayesian image analysis discussed in [40] and [3]. An image $\mathbf{x}$ is a vector of size 256 corresponding to values on the $16 \times 16$ lattice of pixels. Roberts and Sahu [40] model $\mathbf{x}$ using a Gaussian prior density $g$ given by

$$g(\mathbf{x}) \propto e^{-\delta \sum_{i \sim j}(x_i - x_j)^2},$$

where $\delta$ is a constant and $i \sim j$ if $x_i$ and $x_j$ are neighbors. Suppose we observe a corrupted image **y** instead of **x**. The value $y_i$ for each pixel $i$ follows an independent Gaussian density with mean $x_i$ and variance $\sigma^2$. Let $n_i$ denote the number of neighbors of vertex $i$, $1 \leq i \leq 256$. The posterior density of **x** is Gaussian with inverse covariance matrix $Q$ given by

$$Q_{ij} = \begin{cases} 2\delta n_i + \dfrac{1}{\sigma^2}, & i = j, \\ -2\delta, & i \sim j, \\ 0, & \text{otherwise.} \end{cases}$$

Suppose we use the following reversible version of the Gibbs sampling Markov chain to sample from the posterior $\mathcal{N}(\mathbf{0}, Q^{-1})$ density. At every iteration, sample $X_1$ given the other coordinates, then sample $X_2$ given the other coordinates, ..., then sample $X_d$ given the other coordinates twice, then sample $X_{d-1}$ given the other coordinates, ..., and finally sample $X_1$ given the other coordinates. This version of the Gibbs sampler can be expressed in the form (4.8) with

$$A = WL^T,$$
$$V = WDW^T + (D + L^T)^{-1}D(D + L)^{-1},$$

where $D$ and $L$ are the diagonal and lower triangular parts of the matrix $Q$, respectively, and $W = (D + L^T)^{-1}L(D + L)^{-1}$. $A$ satisfies the reversibility condition $AQ^{-1} = Q^{-1}A^T$. For a concrete example, the largest eigenvalue of $A$ when $\delta = 100$ and $\sigma = 0.5$ is 0.9795. Proposition 4.30 tells us that $\chi^2_{\mathbf{0}}(\ell) \leq 10e^{-c}$ for $l = 8.3607 + 12.0620c$ (for any $c \geq 4.8520$), and $\chi^2_{\mathbf{0}}(\ell) \geq e^c/4$ for $l = 8.3607 - 12.0620c$ (for any $c \geq 0$). Figure 9 shows the decrease of the chi-square distance for this chain starting at **0**. Note that eight steps of the Gibbs sampler correspond to $8 \times 2 \times 256 = 4096$ mini sampling steps.

**5. Discussion.** So far, probabilists have come up with various techniques of finding rates of convergence of Markov chains, which can be roughly grouped under five headings:

(a) using the spectral decomposition of a Markov chain,
(b) using Harris recurrence techniques (see [33]),
(c) using probabilistic techniques such as coupling (see [42]), iterated random functions (see [9]) and strong stationary times (see [1, 7]),
(d) using Nash inequalities (see [12]) or logarithmic Sobolev inequalities (see [13]),
(e) using geometric techniques like Poincaré and Cheeger's inequalities (see [8, 30]).

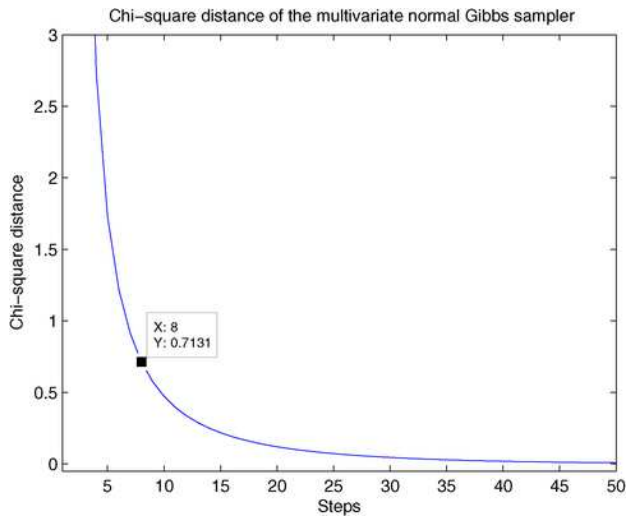

Fig. 9. *Chi-square distance of the multivariate normal Gibbs sampler with $d = 256$, $\delta = 100$, $\sigma = 0.5$, starting from $\mathbf{0}$.*

We use technique (a), that is, using spectral decomposition of a Markov chain, for analyzing all the examples in this paper. An advantage of this technique over the others is that, when applicable, it gives sharp and accurate results. But we require knowledge of all the eigenvalues and eigenfunctions of the Markov chain to apply this technique. In addition, the eigenfunctions have to be suitable for certain algebraic manipulations. This narrows the scope of this technique compared to others. In our examples, the eigenfunctions turn out to be polynomials. We deployed some machinery from the rich field of orthogonal polynomials and were able to compute exact rates of convergence for every Markov chain that we analyzed. But still our success was restricted to classes of natural but special starting points for every Markov chain. We could not appropriately manipulate the distance to stationarity of these Markov chains from a general starting point. Also, exact analysis of the nonreversible case for the multivariate normal autoregressive process remains open. This provides lots of future directions to go, but it is sobering to learn that even in these examples where we know all the eigenvalues and eigenvectors, it is hard, if not impossible, to have an exact analysis from a general starting point. To conclude, the theory of rates of convergence of Markov chains has a long way to go, but it is nice to see standard examples where exact analysis is available by present techniques.

**Acknowledgments.** We thank Persi Diaconis for his valuable guidance and encouragement and Bob Griffiths for help with the kernel polynomials for various distributions covered in this paper. We also thank the referee for a careful review and useful suggestions.

Department of Statistics  
Stanford University  
390 Serra Mall  
Stanford, California 94305  
USA  
E-mail: kdkhare@stanford.edu

Department of Human Genetics  
David Geffen School of Medicine  
University of California, Los Angeles  
Los Angeles, California 90095-1766  
USA  
E-mail: huazhou@ucla.edu